*Статья приурочена к 80-и летию*

*Бориса Теодоровича Поляка*

**© 2016 г. А.В. ГАСНИКОВ, к.ф.-м.н.**

**(ПреМоЛаб МФТИ, ИППИ РАН)**

**П.Е. ДВУРЕЧЕНСКИЙ, к.ф.-м.н.**

**(WIAS Berlin, ИППИ РАН)**

**Ю.Е. НЕСТЕРОВ, д.ф.-м.н.**

**(НИУ ВШЭ, CORE UCL Louvain la Neuve, Belgium)**

# Стохастические градиентные методы с неточным оракулом

В работе предпринята попытка описать современное состояние методов проекции градиента (в том числе прямых методов и методов покомпонентного спуска) решения задач выпуклой стохастической оптимизации с неточным оракулом (неточность неслучайной природы), выдающим стохастический субградиент. Заметная часть приведенных в статье результатов была получена относительно недавно. Цель данной работы – собрать все вместе и посмотреть на разнообразные факты из этой области с единой позиции.

## 1. Введение

В 1960-е годы численные методы выпуклой оптимизации переживали свою первую большую революцию. В работах того времени четко и последовательно развивалась линия градиентных методов. Основополагающим здесь можно признать вклад Бориса Теодоровича Поляка [1, 2], с работ которого во многом и началось активное и повсеместное использование градиентных методов в Советском Союзе. Следующая революция началась в конце 1970-х годов после фундаментальных работ А.С. Немировского, Д.Б. Юдина, Л.Г. Хачияна, N. Karmarkar'a и др. [3, 4]. В монографии [4] была предложена классификация задач выпуклой (и не только) оптимизации по степени гладкости и выпуклости. Были получены нижние оценки для соответствующих классов задач оптимизации с оракулом, выдающим по запросу градиент или стохастический градиент, его компоненту или просто значение функции в точке. Стало понятно, чего в принципе можно достичь. Стали строиться оптимальные методы, см., например, [5–7]. При этом на задачи стали смотреть более пристально с точки зрения теории сложности. Появилась битовая сложность. Была показана полиномиальная разрешимость задач линейного программирования в битовой сложности [3]. Началась разработка полиномиальных методов внутренней точки для задач выпуклой оптимизации на базе метода Ньютона, которая впоследствии привела к созданию общей теории [8–11] и соответствующего пакета CVX [12], способного решать широкий спектр задач выпуклой оптимизации в пространствах размерности до $n \sim 10^4 - 10^5$. Однако, вызовы нового тысячелетия заставляют снова вернуться к градиентным методам. Задачи, которые стали возникать в последние десять лет, отличаются огромными размерностями $n \sim 10^6 - 10^9$. Такие задачи (классифицируемые как задачи large-scale и huge-scale оптимизации) приходят из анализа данных, поиска равновесий в различных сетевых моделях (связанных с компьютерными и транспортными сетями), биоинформатики и многих других областей. Для таких размерностей шаг (итерация) метода Ньютона,[1] становится слишком дорогим, поэтому приходится снова возвращаться к более медленным (в смысле скорости сходимости), но более дешевым (в смысле стоимости одной итерации) градиентным методам (см. [13]). Но для указанных размерностей даже градиентные методы могут испытывать проблемы. В этой связи оказалась очень полезной концепция "заглядывания в черный ящик", т.е. использование структуры задачи с целью ускорения вычислений [14], и использование вместо градиента его легко вычислимой (стохастической) аппроксимации [11]. Как следствие, принято стало считать, что правильный способ эффективно решать ту или иную задачу – это отказаться от общих методов, оптимальных на больших классах, и погружаться в специфику конкретной задачи в надежде ускориться и получить оценки лучше, чем нижние границы [4]. Можно сказать, что началась новая революция. Поток работ на эту тему в основных профилирующих журналах (например, Math. Program.) резко возрос (см., например, обзор [11]). Тем не менее, параллельно стали появляться работы (в том числе работы авторов статьи), показывающие, что многие эффективные методы решения современных задач выпуклой оптимизации в пространствах огромных размеров получаются сочетанием небольшого количества приемов и идей. Цель настоящей работы состоит в том, чтобы собрать воедино набор основных таких идей и показать их связь с некоторыми концепциями 1960-х годов, многие из которых восходят к Б.Т. Поляку. Мы сосредоточимся на оценках числа итераций, требующихся различным методам для решения задачи выпуклой оптимизации с заданной точностью по функции. Эта информация не в полной мере характеризует эффективность метода, но она необходима для последующего его полного исследования. Мы также ограничимся рассмотрением методов проекции градиента [15, 16], в которые, например, не входят очень популярные в последнее время методы условного градиента [11, 17, 18]. В качестве основного инструментария для получения эффективных методов используется метод оценивающих последовательностей, восходящий к работам одного из авторов статьи [7, 10, 14]. Здесь имеются и альтернативные подходы, например, [19–21]. Из-за ограничений на объем статьи и большого количества технических деталей мы ограничимся здесь лишь изложением общей картины. В частности, в статье не приводится псевдокод соответствующих алгоритмов, но, как правило, указываются источники, в которых его можно найти. Мы также не претендуем здесь на полный обзор современного состояния исследований, посвященных градиентным методам. Более того, при ссылках на литературу мы далеко не всегда ссылались на первоисточники, иногда предпочитая ссылаться на удачно написанный более доступный и более современный обзор или монографию.

[1] Заметим, что рассматривать методы более высокого порядка, чем метод Ньютона практически в любой ситуации не имеет смысла, поскольку зависимость числа итераций от желаемой точности решения задачи невозможно сделать лучше, чем у метода Ньютона (в окрестности его квадратичной сходимости), сколько бы старших (равномерно ограниченных) производных функционала не использовалось в методе [4].

## 2. Стохастическая оптимизация

Рассматривается задача выпуклой стохастической оптимизации [2, 22, 23]:

(1) $$f(x) = E_{\xi}\left[ f(x,\xi) \right] \to \min_{x\in Q},$$

где $f(x)$ – выпуклая по $x \in \mathbb{R}^n$ $(n \gg 1)$ функция. Будем называть $\nabla f(x,\xi)$ стохастическим субградиентом функции $f(x,\xi)$ в точке $x$ по первой переменной [24]. Будем считать, что[2] п.н. $\left\|\nabla f(x,\xi)\right\|_2 \le M$, $\nabla = \nabla_x$ и $E_{\xi}$ – перестановочны.[3] Предположим, что $Q$ – выпуклое замкнутое ограниченное множество. Обозначим через $R$ – диаметр множества $Q$: $R = \max_{x,y\in Q}\left\|x-y\right\|_2$. В действительности, достаточно считать, что $R$ – расстояние от точки старта до решения (ближайшего, если решение не единственно) задачи (1) (см. замечание 1). При этом множество $Q$ может быть не ограничено [27].[4] Мы

[2] В действительности [25], здесь и практически в любом другом контексте, где возникает такого типа условия, достаточно требовать, что п.н. выполнено неравенство $\left\|\nabla f(y,\xi) - \nabla f(x,\xi)\right\|_2 \le M$. Это позволяет в ряде случаев понизить оценку константы $M$ и как следствие (см. (2)), ускорить метод. Отметим, что под $\nabla f(x,\xi)$ (аналогично под $\nabla f(y,\xi)$) понимается любой элемент соответствующего стохастического субградиента [24].

[3] Для задач онлайн оптимизации условие перестановочности необходимо записывать в более общем (мартингальном) виде [26].

[4] Впрочем, в случае неограниченного множества $Q$ даже когда выпуклая функция $f(x)$ имеет ограниченную вариацию и равномерно ограниченную норму субградиента на $Q$ (рассматривается детерминированная постановка), мы не имеем никаких гарантий на скорость сходимости используемого метода (под любой метод можно подобрать такую функцию из описываемого класса, что сходимость будет сколь угодно медленной), поскольку не можем никак априорно ограничить расстояние $R$ [4]. Интересный нюанс для выпуклой (но не обязательно сильно выпуклой и гладкой) функции $f(x)$ имеет место, если $f(x)$ задана на ограниченном множестве $Q$ (рассматривается детерминированная постановка). В этом случае размер $R$ множества $Q$ может не входить в оценку необходимого числа итераций (входит вариация

будем считать, что множество $Q$ простой структуры, т.е. на него можно эффективно проектироваться. В работах [27–31] рассматривались различные варианты методов проекции градиента с усреднением и длинными шагами[5] применительно к решению задачи (1). Общая оценка скорости сходимости этих методов есть ($\sigma > 0$ – малый доверительный уровень, $N$ – число итераций метода, на каждой итерации мы можем один раз обратиться к оракулу за субградиентом)[6]

$$P_{x_N}\left( f(x_N) - \min_{x\in Q} f(x) \ge CMR\sqrt{\frac{1+\ln(\sigma^{-1})}{N}} \right) =$$

$$= P_{x_N}\left( E_\xi\left[ f(x_N,\xi)\right] - \min_{x\in Q} E_\xi\left[ f(x,\xi)\right] \ge CMR\sqrt{\frac{1+\ln(\sigma^{-1})}{N}} \right) \le \sigma,$$

где $C$ – константа (здесь и далее константы в основном будут в диапазоне $\sim 10^0 - 10^2$), а случайный вектор $x_N$ – то, что выдает алгоритм (например, метод зеркального спуска [30] или метод двойственных усреднений [27] – сравнительный анализ и описание “физики” этих методов в детерминированном случае проводится в работе [21]) после $N$ итераций. Мы будем называть $x_N$ – $(\varepsilon,\sigma)$-решением задачи (1), если

$$P_{x_N}\left( f(x_N) - \min_{x\in Q} f(x) \ge \varepsilon \right) \le \sigma.$$

Таким образом, для достижения точности по функции $\varepsilon$ и доверительного уровня $\sigma$ методу потребуется (здесь и далее мы будем использовать $\mathrm{O}(\,\cdot\,)$, однако все эти формулы могут быть переписаны с точными константами, что важно, поскольку во многих ситуациях такие оценки используются для формирования критерия останова метода)

(2) $$\mathrm{O}\left( M^2R^2 \ln(\sigma^{-1}) / \varepsilon^2 \right)$$

итераций. На каждой итерации вычисляется стохастический субградиент и осуществляется проектирование.

Отметим, что если использовать метод Монте-Карло, заключающийся в замене исходной задачи (1) следующей задачей

---

функции на этом множестве). Например, это имеет место для метода центра тяжести [2, 4, 11] и для задач, решаемых с относительной точностью [14].

[5] Б.Т. Поляком было показано [28], что такое сочетание позволяет получать эффективные методы для данного класса задач.

[6] Эта оценка неулучшаема с точностью до мультипликативной константы $C$ (при $N \le n$ оценка неулучшаема и в детерминированном случае $f(x,\xi) \equiv f(x)$), см. [4].

(3) $$\frac{1}{N}\sum_{k=1}^{N} f\left(x,\xi_k\right)\to\min_{x\in Q},$$

где с.в. $\xi_k$ – i.i.d., и распределены также как и $\xi$, то для того, чтобы гарантировать, что абсолютно точное решение этой новой задачи является $(\varepsilon,\sigma)$-решением исходной задачи потребуется взять $N$ порядка [24]

$$\mathrm{O}\left(M^2R^2\left(n\ln\left(MR/\varepsilon\right)+\ln\left(\sigma^{-1}\right)\right)/\varepsilon^2\right).$$

Эта наблюдение хорошо поясняет, что подход, связанный с усреднением случайности за счет самого метода более предпочтителен, чем замена задачи (1) ее стохастической аппроксимацией (3).[7] Более предпочтителен не только тем, что допускает адаптивность постановки и легко переносится на онлайн модификации исходной задачи, но, прежде всего, лучшей приспособленностью к большим размерностям.

Здесь важно подчеркнуть фундаментальную идею[8], которую можно усмотреть, например, в [2] и в цикле работ Б.Т. Поляка с Я.З. Цыпкиным [33], о том, что для получения (агрегирования) хороших оценок неизвестных параметров (особенно когда размерность пространства параметров велика) имеет смысл рассматривать задачу поиска оптимальных значений параметра, как задачу стохастической оптимизации и рассматривать выборку как источник стохастических градиентов. Например, истинное значение неизвестного вектора параметров в предположении верности исходной параметрической гипотезы может быть записано как решение задачи стохастической оптимизации [34, 35] (метод наибольшего правдоподобия Фишера)

$$\theta^* = \arg\max_{\theta\in Q} E_\xi\left[L\left(\theta,\xi\right)\right],$$

где $L(\theta,\xi)$ – логарифм функции правдоподобия. Однако решать эту задачу обычными методами мы не можем, потому что математическое ожидание берется по с.в. $\xi$, распределение которой задается неизвестным параметром $\theta^*$. Обойти эту сложность можно, если решать ту же самую задачу

$$E_\xi\left[-L(\theta,\xi)\right]\to\min_{\theta\in Q}$$

методами стохастической оптимизации, получая на каждом шаге новую реализацию (элемент выборки) $\xi_k$ и рассчитывая значения стохастического градиента $\partial L(\theta,\xi_k)/\partial\theta$. То, что выдает алгоритм и будет оценкой вектора неизвестных параметров $\theta^*$. Как правило, дополнительно известно, что $L(\theta,\xi)$ – гладкая и $\mu$-сильно вогнутая (равномерно по $\xi$) функция от $\theta$. Последнее обстоятельство позволяет получить лучшую оценку скорости сходимости по функции [36–39] (в [36, 37] используются специальная модификация метода проекции градиента с усреднением и выбором шагов $h_k = 2\left(\mu\cdot(k+1)\right)^{-1}$ и $h_k = (\mu k)^{-1}$, где $k$ – номер итерации, о подходе [38] и близком к нему подходе [39] будет немного написано в п. 3)

$$\text{(4)}\qquad \mathrm{O}\left(M^2\ln\left(\ln(N)/\sigma\right)/(\mu N)\right),$$

т.е. ($x=\theta$, $f=-L$ , $\bar{C}$ – некоторая константа)

$$P_{x_N}\left(f(x_N)-\min_{x\in Q}f(x)\ge\bar{C}M^2\frac{\ln\left(\ln(N)/\sigma\right)}{\mu N}\right)\le\sigma.$$

---

[7] Особенно ярко это проявляется в случае, бесконечномерных пространств, возникающих в статистической теории обучения (СТО = SLT, Statistical Learning Theory) [32]. Попытка обучиться за счет минимизации эмпирического риска (а именно так можно расшифровать формулу (3) в СТО) может не дать состоятельной оценки/решающего правила, в то время как соответствующий стохастический зеркальный спуск дает состоятельную оценку. Отметим, что в работе [32] приводится достаточно интересный общий результат: в задачах обучения (в частности, в задачах СТО, математической статистики и онлайн обучения) способ получения оптимальных (с точностью до логарифмических факторов) оценок/решающих правил (или, другими словами, способ наискорейшего обучения) базируется на применении соответствующего метода зеркального спуска. Правда, найти “соответствующий метод”, в свою очередь, представляет собой непростую задачу.

[8] Распространяемую и на непараметрическую статистику. Отметим, что начиная с 1980-х годов XX века в этом направлении был цикл работ А.С. Немировского, Б.Т. Поляка и А.Б. Цыбакова, оказавших заметное влияние и на текущие исследования в этой области.

Из неравенства Рао–Крамера [35] ($Q=\mathbb{R}^n$) будет следовать, что оценка (4) – не улучшаемая (с точностью до слагаемого $\ln\left(\ln(N)\right)$). Правда, тут возникают некоторые тонкости, когда мы говорим о неулучшаемости оценок с учетом вероятностей больших отклонений. Строго говоря, классические результаты типа Рао–Крамера, Ван-Трисса и т.п. (см., например, [35]) позволяют лишь говорить о неулучшаемости в смысле сходимости полных математических ожиданий (без вероятностей больших отклонений), и именно в таком смысле можно получить (с помощью методов [11, 37–39]) неулучшаемую (с точностью до мультипликативной константы) оценку:

$$E_{\xi,x_N}\left[f(x_N,\xi)\right]-\min_{x\in Q}E_\xi\left[f(x,\xi)\right]\le\frac{\breve{C}M^2}{\mu N},$$

где $\breve{C}$ – некоторая константа.

Можно обобщить рассмотренную постановку задачи (1) на случай, когда $\left\|\nabla f(x,\xi)\right\|_2$ имеет субгауссовский хвост (определение см., например, в [30]). Оценки (2), (4) при этом сохранят прежний вид. Если же $\left\|\nabla f(x,\xi)\right\|_2^2$ имеет степенной хвост [40], т.е.

$$P\left(\frac{\left\|\nabla f(x,\xi)\right\|_2^2}{M^2}\geq t\right)=\mathrm{O}\left(\frac{1}{t^{\alpha}}\right),$$

где $\alpha > 2$, то[9]

$$P_{x_N}\left(f(x_N)-\min_{x\in Q} f(x)\geq \frac{C_{\alpha}MR}{\sqrt{N}}\left(\sqrt{\ln\left(\sigma^{-1}\right)}+\frac{\left(N/\sigma\right)^{1/\alpha}}{N}\right)\right)\leq\sigma .$$

Если дополнительно $f(x)=E_{\xi}\left[f(x,\xi)\right]$ – $\mu$-сильно выпуклая функция, то (при $\alpha > 2$)

$$P_{x_N}\left(f(x_N)-\min_{x\in Q} f(x)\geq \bar{C}_{\alpha}\frac{M^2}{\mu N}\left(\ln\left(\ln\left(N\ln\left(\sigma^{-1}\right)\right)/\sigma\right)+\frac{\left(\ln\left(N/\sigma\right)/\sigma\right)^{2/\alpha}}{N^{2(1-1/\alpha)}}\right)\right)\leq\sigma .$$

Если ничего не известно о $\left\|\nabla f(x,\xi)\right\|_2^2$, кроме неравенства $E_{\xi}\left[\left\|\nabla f(x,\xi)\right\|_2^2\right]\leq M^2$, то по неравенству Маркова

$$P_{x_N}\left(f(x_N)-\min_{x\in Q} f(x)\geq \frac{\hat{C}MR}{\sigma\sqrt{N}}\right)\leq\sigma ,$$

$$P_{x_N}\left(f(x_N)-\min_{x\in Q} f(x)\geq \frac{\breve{C}M^2}{\sigma\mu N}\right)\leq\sigma ,$$

второе неравенство подразумевает $\mu$-сильную выпуклость $f(x)$.

Можно задать вопрос: насколько вообще уместно рассматривать постановки, в которых возникают тяжелые хвосты. Ведь, если мы можем эффективно вычислять значения функции $f(x)=E_{\xi}\left[f(x,\xi)\right]$ в задаче (1), то ни о каких тяжелых хвостах можно не заботиться. Поскольку, выбрав число шагов так, чтобы метод находил $\varepsilon$-решение с вероятностью $\geq 1/2$, запустив $\log_2\left(\sigma^{-1}\right)$ реализаций такого метода и выбрав реализацию с минимальным значением функции в конечной точке алгоритма, мы за дополнительную $\log_2\left(\sigma^{-1}\right)$ плату (мультипликативную) получим с вероятностью $1-\sigma$ среди выданных ответов хотя бы одно $\varepsilon$-решение [31, 39]. Однако предположение о возможности эффективно вычислять значения функции (при условии трудной вычислимости ее градиента), как правило, не встречается на практике. В некотором смысле типичным тут является пример 1 (см. ниже) вычисления вектора PageRank (при $n \sim 10^9$). Собственно, искусственность ситуации, в которой значение функции легко вычислимо, а градиент нет, неплохо соответствует философии быстрого автоматического дифференцирования (БАД) [41, 42]. Согласно теории БАД, если мы можем посчитать значение функции, то мы можем не более чем в 4 раза дороже посчитать и ее градиент.[10] Как следствие, если мы можем эффективно вычислить значение $f(x)$, то, как правило, мы и $\nabla f(x)$ можем эффективно вычислить. Тогда и на исходную задачу (1) можно смотреть уже не как на задачу стохастической оптимизации, а как на обычную задачу выпуклой оптимизации, что может существенно ускорить ее решение (см. п. 3 ниже). Впрочем, во многих интересных приложениях отмеченный прием (амплификация), как правило, весьма успешно работает [43, 44], поскольку время работы метода, как правило, оказывается заметно большим, чем расчет значения функции.

[9] Приводимые ниже неравенства стоит понимать так, что $x_N$ выдается методом [27, 30], а в сильно выпуклом случае, методом [37–39]. При этом для оценок вероятностей больших уклонений в случае тяжелых хвостов требуется некоторые оговорки и уточнения. К сожалению, мы не смогли найти соответствующий выписанным оценкам (в случае тяжелых хвостов) источник литературы. Приведенные здесь нами формулы нуждаются в дополнительной проверке.

Отметим также (следуя А.С. Немировскому), что с помощью концепции неточного оракула (см. п. 3 ниже) мы можем редуцировать задачу с тяжелыми хвостами $\|\nabla f(x,\xi)\|_2^2$ и компактным множеством $Q$ к ситуации, когда п.н. $\|\nabla f(x,\xi)\|_2 \le M(\varepsilon)$. Для этого нужно "обрезать" стохастический градиент

$$\nabla f(x,\xi) := \begin{cases} \nabla f(x,\xi), & \|\nabla f(x,\xi)\|_2 \le M(\varepsilon) \\ M(\varepsilon)\dfrac{\nabla f(x,\xi)}{\|\nabla f(x,\xi)\|_2}, & \|\nabla f(x,\xi)\|_2 > M(\varepsilon) \end{cases}.$$

Константа $M(\varepsilon)$ подбирается оптимальным образом, исходя из желаемой точности $\varepsilon$. Чем больше $M(\varepsilon)$, тем меньше смещение (bias) обрезанного стохастического градиента, как следствие, тем точнее можно восстановить решение исходной задачи, но при этом возрастает необходимое число итераций (см. (2), (4), в которые входит константа $M = M(\varepsilon)$). Оптимальный выбор этой константы (с точностью до логарифмического фактора) дает приведенные выше оценки.

Все сказанное выше[11] обобщается и на другие прокс-структуры [4] (не обязательно евклидовы, когда выбирается прокс-функция $d(x) = \|x\|_2^2/2$), согласно которым

---

[10] Это легко понять в случае $f(x) = \langle c, x\rangle$. В случае, когда $f(x)$ – многочлен, это также несложно понять (Баур–Штрассен). В общем случае рассуждения аналогичны.

[11] В сильно выпуклом случае (если в прямом пространстве выбрана $q$-норма ($l_q^n$) и прокс-функция $d(x) \ge 0$, $d(x_0) = 0$) в оценку (4) дополнительно входит фактор $\omega = \sup_{x\in Q} 2V(x,x_0)\big/\left(\alpha\|x - x_0\|_q^2\right) \ge 1$ , где $V(x,x_0)$ определяется через $d(x)$ в замечании 1 (при $1 \le q \le 2$ удается найти такую прокс-функцию, что

осуществляется (как правило, по явным формулам[12]) “проектирование” на $Q$. Например, для множества $Q = S_n(1)$ ($S_n(R) = \left\{ x \in \mathbb{R}^n : x_i \ge 0,\ i = 1,\dots,n,\ \sum_{i=1}^{n} x_i = R \right\}$ – единичный симплекс в $n$-мерном пространстве) часто рассматривается (см. пример 1 вычисления вектора PageRank ниже) KL-прокс-структура: $d(x) = \ln n + \sum_{i=1}^{n} x_i \ln x_i$. Эта прокс-функция $d(x) \ge 0$ сильно выпукла в 1-норме с константой сильной выпуклости $\alpha = 1$ на $S_n(1)$ – в силу неравенства Пинскера [27, 30]. Она “наилучшим” образом подходит для симплекса (с некоторыми оговорками [26, 45]). Выгода от ее использования в том, что норма стохастического субградиента всегда оценивается в сопряженном пространстве к пространству, в котором прокс-функция 1-сильно выпукла. В рассматриваемом случае получается $\left\| \nabla f(x,\xi) \right\|_\infty \le M$, что в типичных ситуациях дает оценку константы $M$ в $\sim \sqrt{n}$ раз лучше, чем в 2-норме, а плата за это – увеличение оценки размера области (в этой ситуации в оценке числа итераций нужно использовать $R^2 = \max_{x \in Q} d(x)/\alpha$) в $\sim \ln n$ раз. Детали имеются, например, в статье [30]. Интересным также представляется выбор прокс-функции для прямого произведения симплексов [46]. Здесь мы отметим (следуя А.С. Немировскому), что в общем случае оптимальный выбор прокс-структуры (с точностью до умножения на степень логарифма размерности пространства) связан с симметризацией множества $Q$. Выпуклое центрально симметричное множество $B = (Q - Q)/2$ порождает по теореме Колмогорова норму, в которой $B$ является единичным шаром. Далее ищется оптимальная прокс-функция, согласованная с этой нормой. Говоря более формально, ищется такая сильно выпуклая в этой норме функция $d(x) \ge 0$ с константой сильной выпуклости $\alpha \ge 1$, чтобы число $R^2 = \max_{x \in Q} d(x)/\alpha$ было минимально возможным. Если $Q = B_2^n(1)$ – единичный евклидов шар, то значение $R^2 \le 1$, т.е. не зависит от размерности пространства $n$, но если $Q = B_\infty^n(1)$ – единичный шар в $l_\infty^n$ норме, то $R^2 = \Omega(n)$ (т.е. существует такое число $\chi$, что при достаточно больших значениях $n$ имеет место неравенство $R^2 \ge \chi n$, причем можно добиться того, что $R^2 = O(n)$). Как будет видно из замечания 2 (на примере когда $Q = B_\infty^n(1)$), выбор $l_\infty^n$ нормы не всегда приводит к

---

$\omega = O(\ln n)$, см. замечание 2), где $\alpha$ – константа сильной выпуклости $d(x)$ на $Q$ в $q$-норме [38]. Отметим, что при этом константы в отношение $M^2/\mu$ в оценке (4), считаются относительно $q$-нормы.

[12] Впрочем, в подавляющем большинстве случаев даже если нет возможности явно решить задачу проектирования, ее можно эффективно решить приближенно [9] (посредством перехода к двойственной задаче малой размерности). Как правило, при таком способе рассуждений необходимо использовать концепцию неточного оракула (см. п. 3), поскольку рассчитать градиент двойственного функционала можно лишь приближенно. Однако все эти выкладки обычно не изменяют по порядку сложность одной итерации метода, основной составляющей которой является расчет (пересчет) градиента или его стохастического аналога. Некоторые тонкости и оговорки тут возникают в случае разреженных постановок задач [43, 44].

оптимальным во всех смыслах оценкам (аналогичные примеры нам встретятся и в следующих двух пунктах).

**Замечание 1.** Стоит обратить внимание на то, что если выбрана евклидова прокс-структура, то $R^2$ – квадрат евклидова диаметра $Q$. При переходе к другой прокс-структуре в оценках числа итераций в качестве $R^2$ фигурирует прокс-диаметр $Q$ ($\operatorname{diam}(Q) = \max_{x \in Q} d(x)$), поделенный на константу сильной выпуклости $\alpha = \alpha(Q)$ прокс-функции, заданной на $Q$, относительно выбранной нормы в прямом пространстве. Скажем, в случае выбора KL-прокс-структуры, 1-нормы в прямом пространстве и $Q = S_n(r)$, имеем

$$R^2 = \operatorname{diam}(S_n(r))/\alpha(S_n(r)) = r \cdot \operatorname{diam}(S_n(1))/\left(\alpha(S_n(1))/r\right) = r^2 \cdot (\ln n)/1 = r^2 \ln n .$$

Для евклидовой прокс-структуры размер $Q = S_n(r)$ равнялся бы $2r^2$. Отсюда можно сделать вывод (верный и в общем случае), что выбор прокс-структуры имеет целью оптимально учесть структуру множества с точки зрения того как в итоговую оценку числа итераций будет входить размерность пространства, в котором происходит оптимизация. При гомотетичном увеличении/уменьшении множества оценки числа итераций будут меняться одинаково, независимо от выбранной прокс-структуры. Отметим также, что в формуле (2) для прокс-структуры, отличной от евклидовой точнее писать не $R^2 \ln\left(\sigma^{-1}\right)$, где $R^2 = r^2 \ln n$ (приводим для KL-прокс-структуры), а $r^2\left(\ln n + \ln\left(\sigma^{-1}\right)\right) = r^2 \ln(n/\sigma)$. В действительности, в оценки скоростей сходимости (в среднем, но не в оценки вероятностей больших уклонений, см. замечание 4) всех упомянутых в данной статье методов (кроме обычного (прямого) градиентного метода и метода Франк–Вульфа) входит не прокс-диаметр множества $Q$, на котором происходит оптимизация (если $Q = \mathbb{R}^n$ прокс-диаметр будет бесконечным), а брэгмановское "расстояние" $V(x_*, x_0)$ от решения $x_*$ до точки старта $x_0$ (часто выбирают $x_0 = \arg\min_{x \in Q} d(x)$, $d(x_0) = 0$, $\nabla d(x_0) = 0$ [27]), где $V(x,y) = d(x) - d(y) - \langle \nabla d(y), x - y \rangle$.

**Замечание 2.** Пусть $Q = B_q^n(1)$ – единичный шар в $q$-норме или, в более общем случае, $Q$ содержится в $B_q^n(1)$. Относительно оптимального выбора нормы и прокс-структуры можно заметить следующее (см., например, [4, 9, 47]): если $q \ge 2$, то в качестве нормы оптимально выбирать $\| \ \|_2$ (2-норму) и евклидову прокс-структуру. Определим $q'$ из $1/q + 1/q' = 1$. Пусть $1 \le q \le 2$, тогда $q' \ge 2$. Если при этом $q' = o(\log n)$, то оптимально выбирать $\| \ \| = \| \ \|_q$, а прокс-структуру задавать прокс-функцией $d(x) = \frac{1}{2(q-1)}\|x\|_q^2$. Во всех этих случаях $R^2 = \mathrm{O}(1)$. Для $q' \ge \Omega(\log n)$, выберем $a = 2\log n/(2\log n - 1)$, $\| \ \| = \| \ \|_a$, а прокс-структуру будем задавать прокс-функцией $d(x) = \frac{1}{2(a-1)}\|x\|_a^2$. В этом случае $R^2 = \mathrm{O}(\log n)$. Не сложно проверить, что для единичного симплекса, вложимого в единичный шар в 1-норме, выбор соответствующих прокс-структур из замечаний 1, 2 приводит к одинаковым оценкам числа итераций в категориях $\mathrm{O}(\ )$. В частности, для случая когда $Q = B_\infty^n(1)$, выбор 2-нормы и евклидовой прокс-структуры приводит к оценке (далее в замечании речь идет только об оценке (2)) 1) $\mathrm{O}\left(M_2^2 n \ln\left(\sigma^{-1}\right)/\varepsilon^2\right)$ вместо 2) $\mathrm{O}\left(M_\infty^2 n \ln\left(\sigma^{-1}\right)/\varepsilon^2\right)$ (здесь $E_\xi\left[\|\nabla f(x,\xi)\|_1^2\right] \le M_\infty^2$), получаемой при выборе $l_\infty^n$ нормы в прямом пространстве. Аналогично вышенаписанному можно отметить, что в типичных ситуациях оценка 2 может быть в $\sim n$ раз хуже оценки 1. Тем не менее, оценка 2 $O\left(M_\infty^2 n\ln\left(\sigma^{-1}\right)/\varepsilon^2\right)$ не улучшаема в общем случае.[13] Потому что в общем случае нет гарантий, что $M_2^2 \ll M_\infty^2$, а если такие гарантии есть, то это уже сужает класс функций, для которого получена нижняя оценка с константой $M_\infty^2$.

Подчеркнем, что приведенные здесь оценки (2), (4) (в детерминированном случае при дополнительном условии, что требуемое число итераций для достижения точности $\varepsilon$ удовлетворяет неравенству $N(\varepsilon)\le n$ [4]) без дополнительных предположений являются неулучшаемыми (с точностью до мультипликативных констант) для класса задач стохастической оптимизации (1) и негладких детерминированных задач. Причем дополнительная гладкость функционала задачи (1) в стохастической постановке в общем случае не приводит к улучшению приведенных оценок (2), (4). Если делать дополнительные предположения о малости случайного шума (low noise conditions), то приведенные оценки можно улучшать (см. п. 3). Один пример того, как можно устанавливать неулучшаемость оценок был рассмотрен выше, следуя [2] (на основе неравенства Рао–Крамера), в общем случае следует смотреть монографию [4] и [47]. Отметим, что в работе [47] показывается, что для задач стохастической оптимизации (1) при оптимизации на шарах в $q$-норме оценки типа (2), даваемые методами зеркального спуска с выбором прокс-структуры согласно замечанию 2, соответствуют с точностью до логарифмического фактора нижним оценкам.

Следует, однако, различать задачи стохастической оптимизации и задачи, в которые мы сами искусственно привносим случайность (используя рандомизацию) с целью уменьшения числа арифметических операций на одну итерацию метода [31, 43, 44, 48]. К последнему можно отнести случай, когда (негладкий) выпуклый функционал в задаче является детерминированным, но представляет собой трудно вычислимый интеграл (сумму), зависящую от (оптимизируемых) параметров, который может быть компактно представлен в виде математического ожидания по некоторой простой вероятностной мере. Тогда выгоднее вычислять на каждой итерации метода стохастический градиент, существенно экономя на вычислениях на каждом шаге и лишь немного теряя на логарифмическом увеличении числа шагов ($\sim \ln\left(\sigma^{-1}\right)$). Подробнее об этом подходе будет сказано ниже в примере 3. Ярким примером на эту тему является Google problem (PageRank). По-видимому, одними из первых на эту задачу посмотрели в указанном выше контексте А.В. Назин и Б.Т. Поляк в работе [49], см. также [44, 50–52].

---

[13] Общий результат здесь такой [4, 47]. Пусть необходимо найти минимум выпуклой функции $f(x)$ на множестве $Q=B_q^n(R)$. Оракул выдает несмещенные стохастические субградиенты со свойством $E_\xi\left[\left\|\nabla f(x,\xi)\right\|_{q'}^2\right]\le M_q^2$ $(1/q+1/q'=1)$. Тогда для того, чтобы найти такую точку $x^N$, что $E\left[f\left(x^N\right)\right]-\min_{x\in Q} f(x)\le\varepsilon$, необходимо обратиться к оракулу не менее $N=c_q\, M_q^2R^2/\varepsilon^{\max(2,q)}$, при $N\ll n$, где $c_q=O(\ln n)$ (эта оценка $c_q$ становится точной при $q\to 1+0$); и не менее $N=c_q\, M_q^2R^2n^{1-2/\max(2,q)}/\varepsilon^2$ раз, при $N\gg n$. В детерминированном случае (когда оракул выдает субградиент $\left\|\nabla f(x)\right\|_{q'}\le M_q$) последняя оценка примет вид $N=cn\ln\left(M_qR/\varepsilon\right)$, при $N\gg n$.

**Пример 1 (PageRank).** Задача поиска вектора PageRank $p$ из уравнения $P^T p = p$ ($P$ – стохастическая матрица по строкам матрица), сводится [51, 52] к негладкой задаче выпуклой оптимизации (седловой задаче)

$$\max_{u \in S_n(1)} \left\langle u, P^T p - p \right\rangle \to \min_{p \in S_n(1)}.$$

Перепишем эту задачу в общем виде

$$\min_{x \in S_n(1)} \max_{y \in S_n(1)} \left\langle y, Ax \right\rangle,$$

где матрица $A$ большого размера $n \times n$ (вообще говоря, неразреженная) с элементами, ограниченными по модулю числом $M = 1$. Ключевое наблюдение для решения этой задачи состоит в том [30, 49], что:

$$Ax = E_{i[x]}\left[ A^{\langle i[x] \rangle} \right],$$

где $A^{\langle i \rangle}$ – $i$-й столбец матрицы $A$, вектор $x \in S_n(1)$, а с.в. $i[x]$ имеет категориальное распределение с вектором параметров $x$. Важным следствием является тот факт, что левая часть равенства, $Ax$, вычисляется за $\mathrm{O}\left(n^2\right)$ арифметических операций, а выражение, стоящее в правой части под математическим ожиданием, $A^{\langle i[x] \rangle}$ – всего лишь за $\mathrm{O}(n)$ арифметических операций. Используя это наблюдение (и аналогичное для умножения матрицы $A$ на вектор-строку слева), можно показать, что (рандомизированный) метод зеркального спуска [30] (с KL-прокс-структурой) и стохастическим градиентом по $x$ равным $A^{\langle i[x] \rangle}$ (аналогично по $y$) после

$$\mathrm{O}\left( \frac{nM^2 \ln\left(n/\sigma\right)}{\varepsilon^2} \right) = \mathrm{O}\left( \frac{n \ln\left(n/\sigma\right)}{\varepsilon^2} \right)$$

элементарных арифметических операций выдает такие $x \in S_n(1)$ и $y \in S_m(1)$, что

$$\max_{\tilde{y} \in S_m(1)} \tilde{y}^T Ax - \min_{\tilde{x} \in S_n(1)} y^T A\tilde{x} \le \varepsilon$$

с вероятностью $\ge 1 - \sigma$.

Аналогичные рассуждения [31, 48] позволяют получить с такими же затратами $\mathrm{O}\left(n \ln\left(n/\sigma\right)\varepsilon^{-2}\right)$ такой вектор $x \in S_n(1)$, что $\|Ax\|_\infty \le \varepsilon$. Кроме того, если дополнительно известно, что матрица $P$ – разрежена, то можно организовать поиск $(\varepsilon,\sigma)$-решения еще эффективнее – рандомизировать при проектировании на симплекс [26, 48, 52]. Тогда вместо фактора $n$ в оценках общего числа операций $O\left(n+s\ln n\ln\left(n/\sigma\right)\varepsilon^{-2}\right)$ будет фигурировать $s$ – “среднее” число элементов матрицы $P$ (по строкам и столбцам) отличных от нуля (к сожалению, численные эксперименты Антона Аникина показали, что это “эффективное среднее” число на практике часто близко к максимальному по строкам и столбцам, т.е. от этого подхода можно получить гарантированную выгоду, только если имеет место равномерная разреженность матрицы по строкам и столбцам [48]).

Отметим, что в определенных ситуациях (например, при условии $n \gg \varepsilon^{-2}$ – типичном для задач huge-scale оптимизации) такому рандомизированному методу потребуется использовать относительно небольшое количество элементов матрицы $A$ за все время работы, в то время как для класса детерминированных алгоритмов потребуется считать как минимум половину элементов матрицы $A$ [3] для $\varepsilon = 0.1$.

Хочется также отметить, что на задаче из примера 1 можно продемонстрировать большую часть современного инструментария, необходимого для решения задач huge-scale оптимизации. Так, в случае разреженной матрицы $A$ для решения поставленной негладкой задачи выпуклой оптимизации (и многих других) хорошо подходит метод Б.Т. Поляка [2, 51], работающий по нижним оценкам (2) (функционал негладкий) и при этом учитывающий разреженность $A$ при пересчете градиента [51]. Другой подход [50] (задача поиска вектора PageRank сводится к минимизации другого функционала), также нашедший широкое применение [43, 53–55], связан с заменой градиентного спуска на покомпонентный спуск. Такая замена увеличивает в среднем число итераций всегда не больше (а, как правило, намного меньше) чем в $n$ раз, но зато (благодаря разреженности) происходит экономия при пересчете одной компоненты градиента, как правило (но не всегда – особенности возникают в разреженных задачах), в $n$ раз по сравнению с расчетом полного градиента. В результате получается выгода, которая при определенных условиях может сократить объем вычислений в $\sim\sqrt{n}$ раз (см., например, [43]). Поясним это следующим примером [48], который можно понимать как вариацию неускоренного варианта покомпонентного метода с выбором максимальной компоненты [50].

**Пример 2 (разреженный PageRank).** Задача поиска вектора PageRank также может быть сведена к следующей задаче выпуклой оптимизации (далее для определенности будем полагать $\gamma = 1$, в действительности, по этому параметру требуется прогонка)

$$f(x) = \frac{1}{2}\|Ax\|_2^2 + \frac{\gamma}{2}\sum_{k=1}^{n}\left(-x_k\right)_+^2 \to \min_{\langle x,e\rangle=1}$$

где как и в примере 1 $A = P^T - I$, $I$ – единичная матрица, $e = (1,\ldots,1)^T$,

$$(y)_+ = \begin{cases} y, y \ge 0 \\ 0, y < 0 \end{cases}.$$

При этом мы считаем, что в каждом столбце и каждой строке матрицы $P$ не более $s \ll \sqrt{n}$ элементов отлично от нуля ($P$ – разрежена). Эту задачу предлагается решать обычным градиентным методом[14], но не в евклидовой норме, а в 1-норме (см., например, [21]):

$$x_{k+1} = x_k + \underset{h:\langle h,e\rangle=0}{\arg\min}\left\{ f(x_k) + \langle \nabla f(x_k), h\rangle + \frac{L}{2}\|h\|_1^2 \right\},$$

где $L = \max_{i=1,\ldots,n}\left\|A^{\langle i\rangle}\right\|_2^2 + \gamma \le 3$ ($A^{\langle i\rangle}$ – $i$-й столбец матрицы $A$). Для достижения точности $\varepsilon^2$ по функции потребуется сделать $\mathrm{O}\left(LR^2/\varepsilon^2\right) = \mathrm{O}\left(1/\varepsilon^2\right)$ итераций [14]. Не сложно проверить, что пересчет градиента на каждой итерации заключается в умножении $A^T Ah$, что может быть сделано за $\mathrm{O}\left(s^2\right)$. Связано это с тем, что вектор $h$ всегда имеет только две компоненты отличные от нуля (такая разреженность получилась благодаря выбору 1-нормы), причем эти компоненты соответствуют $\underset{i=1,\ldots,n}{\arg\min}\,\partial f(x_k)/\partial x^i$ и $\underset{i=1,\ldots,n}{\arg\max}\,\partial f(x_k)/\partial x^i$, что пересчитывается (при использовании специального двоичного дерева (кучи) для поддержания максимальной и минимальной компоненты градиента [51]) за $\mathrm{O}\left(s^2 \ln n\right)$ (логарифмический фактор можно ослабить, если использовать, например, фибоначчиевы или бродалевы кучи [48]). Таким образом, общая трудоемкость предложенного метода будет $\mathrm{O}\left(n + s^2 \ln n/\varepsilon^2\right)$, что заметно лучше многих известных методов [52]. Стоит также отметить, что функционал, выбранный в этом примере, обеспечивает намного лучшую оценку $\|Ax\|_2 \le \varepsilon$ по сравнению с функционалом из примера 1, который (в варианте [31]) обеспечивает $\|Ax\|_\infty \le \varepsilon$. Наилучшая (в разреженном случае без, условий на спектральную щель матрицы $P$ [52]) из известных нам на данный момент оценок $\mathrm{O}\left(s \ln n \ln(n/\sigma)/\varepsilon^2\right)$ [26, 52] для $\|Ax\|_\infty$ может быть улучшена приведенной в этом примере оценкой, поскольку, как уже отмечалось ранее, $\|Ax\|_2$ может быть (и так часто бывает) в $\sim\sqrt{n}$ раз больше $\|Ax\|_\infty$, а $s \ll \sqrt{n}$.

[14] Выписанная далее оценка скорости сходимости (на число итераций) – неулучшаема с точностью до мултипликативного фактора. Речь идет не об оптимальности метода на классе гладких задач на симплексе, а о том, что конкретно для этого метода такая оценка если и может быть улучшена, то лишь на мультипликативный фактор. Это замечание касается практически всех известных сейчас градиентных методов. Показывается это приблизительно также (даже еще проще), как и в случае оптимальности оценок на классах [4]: строится конкретные примеры (семейства) функций.

Заметим, что в решении могут быть маленькие отрицательные компоненты. Также численные эксперименты показали [48], что для достижения выписанных оценок требуется препроцессинг (в нашем случае он заключается в представлении матрицы по строкам в виде списка смежности: в каждой строке отличный от нуля элемент хранит ссылку на следующий отличный от нуля элемент, аналогичное представление матрицы делается и по столбцам). Заметим, что препроцессинг помогает ускорять решение задач не только в связи с более полным учетом разреженности постановки, но и, например, в связи с более эффективной организацией рандомизации [31, 43, 50].

Пример 2 также характерным образом демонстрирует, как используется разреженность (см. также [44, 51, 56]). Обратим внимание на то, что число элементов в матрице $P$, отличных от нуля, даже при наложенном условии разреженности (по строкам и столбцам), все равно может быть достаточно большим $sn$. Удивляет то, что в оценке общей трудоемкости это число не присутствует. Это в перспективе (при правильной организации работы с памятью) позволяет решать задачи огромных размеров. Более того, даже в случае небольшого числа не разреженных ограничений вида $\langle a_i, x\rangle = b_i$, $i = 1,..,m = \mathrm{O}(1)$, можно “раздуть” пространство (не более чем в два раза), в котором происходит оптимизация (во многих методах, которые учитывают разреженность такое раздутие не приведет к серьезным затратам), и переписать эту систему в виде $Ax = b$, где матрица будет иметь размеры $\mathrm{O}(n)\times\mathrm{O}(n)$, но число отличных от нуля элементов в каждой строке и столбце будет $\mathrm{O}(1)$. Таким образом, допускается небольшое число “плотных” ограничений.

Заметим, что если применить метод условного градиента [17] (Франк–Вульфа) к задаче из примера 2, то общая трудоемкость (для точности $\varepsilon^2$, как и в примере 2) будет [48, 56]

$$\mathrm{O}\left(n + \frac{s^2 \ln\left(2 + n/s^2\right)}{\varepsilon^2}\right).$$

В связи со сказанным выше, заметим, что задача может быть не разрежена, но свойство разреженности появляется в решении при использовании метода Франк–Вульфа, что также может заметно сокращать объем вычислений в постановках аналогичных примеру 2, но с матрицами $A$, у которой число столбцов на много порядков больше числа строк (см., например, п. 3.3 [11], [57]).

Приведем еще один пример, подсказывающий, как следует решать задачу (3), полученную из (1) с применением идеи метода Монте-Карло.

**Пример 3 (рандомизация суммы).** Пусть необходимо решить задачу выпуклой оптимизации (или ее композитный вариант, см., например, замечание 6)

(5) $$f(x)=\frac{1}{N}\sum_{k=1}^{N} f_k(x)\to\min_{x\in Q},$$

где $f_k(x)$ – негладкие выпуклые функции с ограниченной числом $M$ нормой субградиента, $Q$ – выпуклое замкнутое множество простой структуры (можем эффективно на него проектироваться, согласно заданной прокс-функции) прокс-диаметра $R$. Введем новую функцию

$$f(x,\xi)=\begin{cases} f_1(x), \text{ с вероятностью } 1/N \\ \ldots\ldots\ldots\ldots\ldots\ldots\ldots\ldots\ldots\ldots \\ f_N(x), \text{ с вероятностью } 1/N \end{cases}.$$

Ее стохастический субградиент легко вычислить. Для этого разыгрывается за $\mathrm{O}(\ln N)$ с.в. $\xi$, принимающая значения $1,\ldots.,N$ с равными вероятностями (см., например, [52]). Затем считается субградиент $f_\xi(x)$ (и выполняется прокс-проектирование на $Q$). Как уже отмечалось ранее, можно найти $(\varepsilon,\sigma)$-решение так понимаемой задачи (5) за

$$\mathrm{O}\left(\frac{M^2R^2\ln(\sigma^{-1})}{\varepsilon^2}\right)$$

итераций, со стоимостью одной итерации равной $\mathrm{O}(\ln N)$ + затраты на вычисления субградиента $f_\xi(x)$ + затраты на вычисление проекции. Если решать задачу без рандомизации, то число итераций будет $\mathrm{O}\left(M^2R^2/\varepsilon^2\right)$, строго говоря, здесь $M$ должно быть немного меньше за счет того, что

$$\max_{x\in Q}\left\|\frac{1}{N}\sum_{k=1}^{N}\nabla f_k(x)\right\|_* \le \max_{\substack{k=1,\ldots,N\\ x\in Q}}\left\|\nabla f_k(x)\right\|_*,$$

но мы считаем, что обе части неравенства одного порядка. Зато шаг итерации будет теперь почти в $N$ раз дороже. И если $N\gg 1$ это может оказаться существенным.

Приведенную постановку можно распространить на случай, когда взвешивание функций не равномерное (тогда первое разыгрывание с.в. $\xi$, имеющей категориальное распределение, или приготовление процедуры рандомизации займет $\mathrm{O}(N)$, а все последующие $\mathrm{O}(\ln N)$) и $f_k(x):=E_{\xi_k}\left[f_k(x,\xi_k)\right]$ с равномерно ограниченными (по $k$, $x$ и $\xi$) нормами стохастических субградиентов. При этом все приведенные оценки числа итераций сохранятся. Причем требование равномерной ограниченности норм стохастических субградиентов можно существенно ослабить за небольшую плату (см. выше).

Если на решение задачи (3) теперь посмотреть в контексте описанной рандомизации с $f_k(x) = f(x, \xi_k)$ (здесь $\xi_k$ – не случайная величина, а полученная в методе Монте-Карло $k$-я по порядку реализация с.в. $\xi$), то "все встанет на свои места" в смысле одинаковости (с точностью до логарифмического фактора) двух подходов к решению задачи (1), описанных в начале пункта.

Описанная рандомизация при вычислении субградиента суммы функций, по-видимому, была одной из первых, которые предлагались в стохастической оптимизации [22]. Однако она популярна и по сей день, например, в связи с приложениями к поиску равновесий в транспортных сетях [58–62] и анализу данных [63–65]. В частности, в [11, 43, 59, 66–71] в предположении, что все функции в (5) гладкие с константой Липшица градиента $L$, предложен специальный рандомизированный метод (на базе описанного выше способа рандомизации суммы), в котором число вычислений градиентов слагаемых[15]

$$\mathrm{O}\left(\left(N + \min\left\{LR^2/\varepsilon, \sqrt{NLR^2/\varepsilon}\right\}\right)\left(\ln\left(\Delta f/\varepsilon\right) + \ln\left(\sigma^{-1}\right)\right)\right),$$

где $\Delta f$ разность значения функции в стартовой точке и в минимуме. Эта оценка с точностью до выражения под логарифмом соответствует нижней оценке в классе детерминированных алгоритмов [70, 72]. Если дополнительно имеется еще и $\mu$-сильная выпуклость $f(x)$, то оценку можно переписать следующим образом

$$\mathrm{O}\left(\left(N + \min\left\{L/\mu, \sqrt{NL/\mu}\right\}\right)\left(\ln\left(\Delta f/\varepsilon\right) + \ln\left(\sigma^{-1}\right)\right)\right).$$

Отметим, что вторая оценка переходит в первую при следующей квадратичной регуляризации. К выпуклому функционалу прибавляется регуляризирующее слагаемое $\mu\|x\|_2^2/2$. В результате функционал становится сильно выпуклым и справедлива вторая оценка на число вычислений градиента. Такая регуляризация изменяет исходную целевую функцию на число не больше $\mu R^2/2$ и чтобы итоговая погрешность по исходной функции была порядка $\varepsilon$ нужно выбирать $\mu \simeq \varepsilon/R^2$, и решать регуляризованную задачу с точностью $\varepsilon/2$. При подстановке этого значения во вторую оценку числа вычислений градиента последняя переходит в первую оценку.

Отметим также, что сначала (см., например, [55, 68]) получается результат о сходимости средних[16]

[15] Строго говоря, имеющиеся сейчас рассуждения для второго аргумента минимума [43, 73] позволяют получить только при дополнительных предположениях о структуре задачи оценку, аналогичную приведенной ниже, и то только в категориях общего числа арифметических операций.

[16] Описанная далее конструкция не зависит от того, изначально имела место сильная выпуклость или мы ее искусственно ввели должной регуляризацией.

$$E\left( f\left(x_N\right) - \min_{x\in Q} f\left(x\right)\right) \le \varepsilon,$$

где

$$N = N\left(\varepsilon\right) = \mathrm{O}\left(\left(N + \min\left\{L/\mu, \sqrt{NL/\mu}\right\}\right)\ln\left(\Delta f/\varepsilon\right)\right),$$

Потом из неравенства Маркова получают оценку больших уклонений

$$P\left( f\left(x_{N(\varepsilon)}\right) - \min_{x\in Q} f\left(x\right) \ge \sigma\right) \le \varepsilon/\sigma,$$

которую переписывают в виде

$$P\left( f\left(x_{N(\varepsilon\sigma)}\right) - \min_{x\in Q} f\left(x\right) \ge \sigma\right) \le \varepsilon,$$

где

$$N\left(\varepsilon\sigma\right) = \mathrm{O}\left(\left(N + \min\left\{L/\mu, \sqrt{NL/\mu}\right\}\right)\left(\ln\left(\Delta f/\varepsilon\right) + \ln\left(\sigma^{-1}\right)\right)\right).$$

Мы привели здесь это наблюдение, потому что оно оказывается полезным и во многих других контекстах, в которых рандомизированный метод сходится со скоростью геометрической прогрессии.

При наличии дополнительной структуры у задачи (5) приведенные оценки можно было получить (и даже немного улучшить, например, учитывая разреженность) исходя из рандомизированных покомпонентных методов (например, ALPHA или APPROX [63] или ACRCD* из замечания 8 [43]) для “двойственной” к (5) задаче [43, 55, 59, 73, 74].[17] Заметим также, что в работе [43] показывается, как можно просто получить часть выписанных оценок с помощью метода, работающего по оценкам (8) (см. ниже).

В книге [2] Б.Т. Поляк отмечает, что если рандомизация осуществляется каким-то специальным образом, например, таким, что[18]

$$(6) \qquad E\left[\left\|\nabla f\left(x,\xi\right)\right\|_*^2\right] \le C_n \left\|\nabla f\left(x\right)\right\|_*^2 + \Delta,$$

[17] Строго говоря, построение двойственной задачи предполагает возможность явного выделения в функционале в виде отдельного слагаемого сильно выпуклого композита – желательно сепарабельного.

[18] Если рассматривать приложения методов стохастической оптимизации к СТО [32], а в правой части неравенства вместо $\left\|\nabla f\left(x\right)\right\|_*^2$ писать $\left\|\nabla f\left(x\right)\right\|_*^{1/2}$, то выписанное неравенство будет соответствовать условиям малого шума Цыбакова–Массара, Бернштейна [75].

где $\Delta \geq 0$ некоторая малая погрешность, и в точке минимума $\nabla f(x) = 0$, то приведенные выше оценки (2), (4) можно существенно улучшить. Примеры будут приведены ниже в п. 4 (см. (22)). В частности, в сильно выпуклом случае можно получить геометрическую скорость сходимости. Важно отметить, что при рандомизации, возникающей в покомпонентных спусках, спусках по направлению и безградиентных методах в гладком случае условие (6) выполняется [2, 50, 76]. Мы вернемся к этому кругу вопросов в п. 4. Описанная же выше конструкция (с довольно грубым неравенством Маркова) используется в данном контексте [50, 76] для (точной!) оценки больших уклонений. Причем за счет регуляризации функционала, о которой было сказано выше, все это переносится и просто на гладкий случай без предположения сильной выпуклости.

## 3. Стохастические градиентные методы с неточным оракулом

В этом пункте мы опишем, что можно получить, если дополнительно известно, что $f(x)$ – гладкая по $x$ функция, с константой Липшица градиента $L$ и(или) сильно выпуклая с константой $\mu \geq 0$, но вычисление стохастического градиента на каждом шаге происходит с неконтролируемой неточностью $\delta$, вообще говоря, не случайной природы.[19]

**Замечание 3.** И гладкости и сильной выпуклости можно добиться искусственно. Как уже отмечалось в п. 2, сильная выпуклость всегда легко получается регуляризацией функционала в исходной задаче. Как правило, это не дает ничего нового с точки зрения выписанных оценок (и даже может ухудшать эти оценки на логарифмический фактор), но в ряде специальных случаев (см. ниже) это может давать определенные преимущества. Кроме того, такая регуляризация иногда просто необходима для корректности постановки. Это связано с тем, что в общем случае даже для гладких детерминированных выпуклых задач мы можем гарантировать сходимость итерационного метода лишь по функции, но не по аргументу. Для сходимости по аргументу нужна сильная выпуклость функционала, которую и обеспечивают должной регуляризацией (см., например, конец п. 2) – при этом сходимость по аргументу имеет место к решению регуляризованной задачи. Идея регуляризации используется в популярном методе двойственного сглаживания [77] (регуляризация двойственной задачи с целью улучшения гладких свойств прямой). В отличие от прямой регуляризации, эта техника хорошо работает только для вполне конкретных задач, имеющих определенную (седловую – Лежандрову) структуру (модель), когда исходная задача имеет явное двойственное представление (см. пример 4), введя в которое регуляризацию, можно явно (эффективно) пересчитать во что превратится исходная прямая задача. Другой пример сглаживания будет приведен в п. 4.

Сформулируем более точно предположения об оракуле, выдающем стохастический градиент, следуя [78, 79].[20]

[19] Особое внимание таким постановкам стали уделять после выхода книг [2, 33]. В них обстоятельно изучается "влияние помех", в том числе не случайной природы, на методы выпуклой оптимизации.

[20] В работе [78] собрано много различных мотиваций такому предположению (определению), обобщающему классическую концепцию $\delta$-субградиента [2]. В определенном смысле это предположение 1 наиболее общее и, одновременно, наиболее точно отражающее спектр всевозможных приложений [43, 55, 60–62].

**Предположение 1.** $(\delta, L, \mu)$*-оракул выдает (на запрос, в котором указывается только одна точка* $x$*) такую пару* $\left(F(x,\xi), G(x,\xi)\right)$ *(с.в.* $\xi$ *независимо разыгрывается из одного и того же распределения, фигурирующего в постановке (1)), что для всех* $x \in Q$ *ограничена дисперсия*

$$E_\xi\left[\left\|G(x,\xi)-E_\xi\left[G(x,\xi)\right]\right\|_*^2\right] \le D,$$

*и для любых* $x, y \in Q$

$$\frac{\mu}{2}\|y-x\|^2 \le E_\xi\left[f(y,\xi)\right]-E_\xi\left[F(x,\xi)\right]-\left\langle E_\xi\left[G(x,\xi)\right], y-x\right\rangle \le \frac{L}{2}\|y-x\|^2+\delta .$$

Из недавних результатов [78–87] можно получить общий метод (мы приводим огрубленный вариант оценки времени работы этого метода для большей наглядности), с такими оценками скорости сходимости[21]

$$\text{(7)} \qquad \min\left\{O\left(\frac{LR^2}{N^{p+1}}+\sqrt{\frac{DR^2}{N}}+N^p\delta\right), O\left(LR^2\exp\left(-\Upsilon N\cdot\left(\frac{\mu}{L}\right)^{\frac{1}{p+1}}\right)+\frac{D}{\mu N}+\left(\frac{L}{\mu}\right)^{\frac{p}{p+1}}\delta\right)\right\},$$

где $\Upsilon \ge 1$ – некоторая константа (можно считать $\Upsilon = O(\ln n)$), а параметр $p \in [0,1]$ подбирается “оптимально” перед запуском метода исходя из масштаба шума $\delta$. Для лучшего понимания оценки (7) полезно ее переписать в еще более огрубленном виде[22]

$$\min\left\{O\left(\frac{LR^2}{N^{p+1}}+\sqrt{\frac{DR^2}{N}}+N^p\delta\right), O\left(LR^2\exp\left(-\Upsilon N\cdot\left(\frac{\mu}{L}\right)^{\frac{1}{p+1}}\right)+\frac{D}{\mu N}+N^p\delta\right)\right\}.$$

---

[21] Оценки характеризуют достигнутую в среднем точность (по оптимизируемому функционалу) после $N$ итераций. При этом в случае когда минимум достигается на втором аргументе (выгодно использовать факт наличия $\mu$-сильной выпуклости) под $N$ правильнее понимать не число итераций, а число обращений к $(\delta, L, \mu)$-оракулу [87]. Отметим, что при $p=0$ оценку можно сделать непрерывной по параметру $\mu \ge 0$ (см. [79]). Также заметим, что в метод (например, в размер шагов) не входит требуемое число итераций (или желаемая точность – одно через другое выражается). Таким образом, можно говорить об адаптивности метода. Отметим, что за это не приходится дополнительно платить логарифмическую плату [27]. Отметим также, что при $p=1$ метод наихудшим (а при $p=0$ наилучшим) способом (среди всех разумных вариаций градиентного метода) накапливает неточность в вычислении градиента. Это переносится и на негладкие задачи (см. далее).

[22] Насколько нам известно, для всех методов, которые используют только градиент и значение функции (или их стохастические аналоги) накопление шума методом со скоростью $N^p\delta$ с $p \in [0,1]$ – является общим местом.

Этот общий метод – есть в некотором смысле “выпуклая комбинация” двойственного градиентного метода (DGM) и быстрого градиентного метода[23] (FGM) [83, 87], оценки скорости сходимости для которых имеют соответственно вид:

$$\text{(DGM)} \qquad \min\left\{ O\left( \frac{LR^2}{N} + \sqrt{\frac{DR^2}{N}} + \delta \right), O\left( LR^2 \exp\left( -\Upsilon_1 N \frac{\mu}{L} \right) + \frac{D}{\mu N} + \delta \right) \right\},$$

$$\text{(FGM)} \qquad \min\left\{ O\left( \frac{LR^2}{N^2} + \sqrt{\frac{DR^2}{N}} + N\delta \right), O\left( LR^2 \exp\left( -\Upsilon_2 N \sqrt{\frac{\mu}{L}} \right) + \frac{D}{\mu N} + \sqrt{\frac{L}{\mu}}\delta \right) \right\}.$$

Комбинируя эти два метода можно непрерывно настраиваться (оптимально подбирая метод, регулируя $p \in [0,1]$) на шум (известного масштаба). В этой связи также полезно отметить (аналогичный факт имеет место и для покомпонентного варианта FGM [43]), что FGM есть специальная выпуклая комбинация прямого градиентного метода (PGM), оценки скорости сходимости которого совпадают с оценками DGM, и метода зеркального спуска / двойственных усреднений [21, 92] (по-видимому, здесь вместо зеркального спуска можно использовать и метод из работы [93]). Нельзя в этой связи не обратить внимание на то, что комбинация двух методов привела к новому методу, работающему лучше, чем каждый из методов по отдельности. Отметим здесь также недавнюю работу [94], в которой предлагается общий способ получения ускоренных (быстрых) методов.

Вся последующая часть п. 3 будет посвящена обсуждению этих результатов и их окрестностей.

Прежде всего, заметим, что дисперсию у первого аргумента минимума в (7) можно уменьшать в $m$ раз, запрашивая на одном шаге реализацию стохастического градиента не один раз, а $m$ раз, и заменяя стохастический градиент средним арифметическим [11, 31, 79] (в случае тяжелых хвостов у стохастических градиентов лучше пользоваться более робастными оценками, например, медианного типа [4]).[24] Это имеет смысл делать, если слагаемое, отвечающее стохастичности, доминирует. Важно, что мы при этом не увеличиваем число итераций, и слагаемое $N^p\delta$ остается прежним. Отметим, что число вызовов оракула при этом увеличивается, но тем не менее, в некоторых ситуациях такой подход может оказаться оправданным. Такая игра используется[25] в способе получения второго аргумента оценки (7). В этой связи оценку (7) правильнее переписать следующим образом (здесь $N(\varepsilon)$ – число обращений к $(\delta, L, \mu)$-оракулу, необходимых для достижения в среднем по функции точности $\varepsilon$, индекс 1 соответствует просто выпуклому, а индекс 2 сильно выпуклому случаю):

$$(8)\quad N_1(\varepsilon) = \max\left\{ \mathrm{O}\left(\left(\frac{LR^2}{\varepsilon}\right)^{\frac{1}{p+1}}\right), \mathrm{O}\left(\frac{DR^2}{\varepsilon^2}\right)\right\},\ N_2(\varepsilon) = \max\left\{ \mathrm{O}\left(\left(\frac{L}{\mu}\right)^{\frac{1}{p+1}} \ln\left(\frac{LR^2}{\varepsilon}\right)\right), \mathrm{O}\left(\frac{D}{\mu\varepsilon}\right)\right\}$$

при (условия на допустимый уровень шума, при котором оценки (8) имеют такой же вид, с точностью до $\mathrm{O}(1)$, как если бы шума не было)

$$(9)\qquad \delta_1(\varepsilon) \le \mathrm{O}\left(\varepsilon \cdot \left(\frac{\varepsilon}{LR^2}\right)^{\frac{p}{p+1}}\right),\ \delta_2(\varepsilon) \le \mathrm{O}\left(\varepsilon \cdot \left(\frac{\mu}{L}\right)^{\frac{p}{p+1}}\right).$$

---

[23] Отметим, что при $D = 0$ не улучшаемые оценки, которые дает метод FGM [11], были установлены Б.Т. Поляком [2] для ряда других многошаговых методов (метод тяжелого шарика, сопряженных градиентов). Отличие в том, что тогда оценки были установлены локально. Все приведенные в данной статье оценки – глобальные, т.е. не требуют оговорок о близости точки старта к решению, для гарантии нужной скорости сходимости. Заметим также, что техника установления локальной сходимости основана, как правило, на первом методе Ляпунова [2, 88], в то время как глобальной – на втором [2, 89]. При этом функцию Ляпунова можно искать по непрерывному аналогу итерационного процесса – системе дифференциальных уравнений [89]. Скажем, для обычного градиентного метода это будет система [1] (Коши, 1847): $dx/dt = -\nabla f(x)$. Скорости сходимости у итерационного процесса и его непрерывного аналога могут отличаться. Скажем, непрерывный аналог метода Ньютона сходится за конечное время. Другой пример – метод зеркального спуска [4]. Недавно появилась работа, посвященная и непрерывному аналогу FGM [90], см. также [91].

[24] Этот прием в западной литературе часто называют “mini-batch” [11].

Как уже отмечалось, выписанные оценки (7) ((8), (9)) характеризуют скорость сходимости в среднем. Они с одной стороны не улучшаемы[26] с точностью до мультипликативной константы (см. п. 2 и [4, 47]), а с другой стороны достигаются. Все это (неулучшаемость оценок) справедливо и при $\delta = 0$ и(или) $D = 0$. При этом в случае $D = 0$, $\mu = 0$ необходимо считать, что требуемое число итераций для достижения точности $\varepsilon$ удовлетворяет неравенству $N(\varepsilon) \le n$ [4], в противном случае оценки улучшаемы – метод центров тяжести [4, 11], с оценкой числа итераций типа $\mathrm{O}(n \ln(B/\varepsilon))$, где $|f(x)| \le B$. В терминах больших отклонений возникают оценки, аналогичные тем, которые были приведены в п. 2, см. [87].

Отмеченные результаты переносятся и на прокс-структуры отличные от евклидовой [87]. При этом рассмотрение какой-либо другой $q$-нормы ($l_q$-нормы) в прямом пространстве ($q \ge 1$), отличной от евклидовой, в сильно выпуклом случае (когда минимум достигается на втором выражении в (7)), как правило, не имеет смысла. Связано это с тем, что квадрат евклидовой асферичности $q$-нормы, который может возникать в оценках числа обусловленности прокс-функции в $q$-норме (это число, в свою очередь, оценивает увеличение числа итераций метода при переходе от евклидовой норме к $q$-норме), больше либо равен 1. Равенство достигается на евклидовой норме. Скажем, для 1-нормы эта асферичность оценивается снизу размерностью пространства [14, 38]. Другими словами, действительно, можно выбирать в сильно выпуклом случае $q$-норму (отличную от евклидовой) и получать оценки на число итераций вида (см. (7), (8) и п. 2)

$$O\left(\left(\frac{L}{\mu}\omega\right)^{\frac{1}{p+1}}\ln\left(\frac{LR^2}{\varepsilon}\right)\right),\ \omega=\sup_{x\in Q}\frac{2V(x,x_0)}{\alpha\|x-x_0\|_q^2},$$

где $L$ и $\mu$ считаются относительно $q$-нормы, а $R^2$ – брэгмановское “расстояние” от точки старта до решения (см. замечание 1) . Однако смысла, как правило, в этом нет, поскольку $\omega\geq 1$, а число обусловленности $\chi=L/\mu$ не меньше чем в случае выбора 2-нормы. Например [38], для функции $\|x\|_2^2=x_1^2+...+x_n^2$ в евклидовой норме число обусловленности $\chi=1$, а в 1-норме $\chi=n$. Тем не менее, выгода от использования не евклидовой прокс-структуры в сильно выпуклом случае может быть, если рассматривать задачи композитной оптимизации, в которых сильная выпуклость приходит от композитного слагаемого (см. замечание 6). Так в приложениях, описанных в работах [55, 59], в качестве композитного слагаемого возникает сильно выпуклая в 1-норме энтропийная функция. Отметим, что энтропию при этом нельзя использовать в качестве прокс-функции. Нужно брать (и это можно сделать, см. замечание 2) другую прокс-функцию, соответствующую 1-норме, которая обеспечивает (по-видимому, оптимально возможную) оценку $\omega=O(\ln n)$.

---

[25] Вместе с идеей рестартов [38, 43, 56, 84–87], распространяющей (ускоряющей) практически любой итерационный метод (желательно с явной оценкой необходимого числа итераций $N(\varepsilon)$ для достижения заданной точности $\varepsilon$) на случай сильно выпуклого функционала. Нетривиально здесь то, что при довольно общих условиях при таком распространении сохраняется (и работает уже в условиях сильной выпуклости) свойство оптимальности исходного метода.

[26] В нижнюю оценку во втором выражении под знаком минимума при экспоненте вместо $LR^2$ входит $\mu R^2$, а константа $\Upsilon = 1$ в (7). Впрочем, получить вместо фактора $LR^2$ фактор $\mu R^2$ можно аккуратно проанализировав оценки, даваемые с помощью техники рестартов (см., например, [21, 87]).

Заметим также, что обычный метод FGM в не стохастическом сильно выпуклом случае для задач безусловной оптимизации, в действительности, дает оценку (следует сравнить с (8)) [10]:

$$O\left(\sqrt{\frac{L}{\mu}}\ln\left(\frac{f(x_0)-f(x_*)}{\varepsilon}\right)\right).$$

Поскольку $\nabla f(x_*)=0$, то $f(x_0)-f(x_*)\leq LR^2/2$. Если рассматривается задача условной оптимизации (на выпуклом множестве $Q\subset\mathbb{R}^n$), то, вообще говоря, $\nabla f(x_*)\neq 0$, следовательно, нельзя утверждать, что $f(x_0)-f(x_*)\leq LR^2/2$. В [79, 87] предлагается обобщение классического FGM для класса гладких сильно выпуклых задач, которое фактически позволяет вместо $f(x_0)-f(x_*)$ писать нижнюю оценку $\mu R^2/2\leq f(x_0)-f(x_*)$ в том числе для задач условной оптимизации. Заметим, что это же наблюдение справедливо для описываемых в данной работе методов (мы не стали писать $\mu R^2$ вместо $LR^2$ в (7) и далее для сохранения непрерывности выписанных оценок по $\mu$, т.е. чтобы делать меньше оговорок о переключениях с сильно выпуклого случая на выпуклый при малых значениях $\mu$).

**Замечание 4.** Отметим, что пока нам не известно (для произвольной прокс-структуры, отличной от евклидовой) строгое обоснование оценок (7) ((8), (9)) с вероятностями больших отклонений для случая не ограниченного множества $Q$. В известном нам способе получения оценок вероятностей больших уклонений (см., например, [79, 87]), к сожалению, явно используется предположение об ограниченности множества $Q$. С другой стороны, для используемых в статье неускоренных методов (кроме Франк–Вульфа и кроме PGM в варианте [21], для PGM в варианте [79] все хорошо) оценки на скорость сходимости обычно получаются в следующем виде [14, 27, 60, 61, 79, 87, 92]:

$$\sum_{k=0}^{N}\lambda_k\cdot\left(f(x_k)-f(x_*)\right)\le V(x_*,x_0)-V(x_*,x_{N+1})+\sum_{k=0}^{N}\lambda_k\left\langle G(x_k,\xi_k)-\nabla f(x_k),x_*-x_k\right\rangle+$$
$$+\tilde{\Delta}_N\left(\{\lambda_k\}_{k=0}^{N},\left\{\left\|G(x_k,\xi_k)-\nabla f(x_k)\right\|_*^2\right\}_{k=0}^{N},\delta\right),\ \{\lambda_k\}\ge 0$$

или, в случае ускоренных методов (к которым относится FGM и его производные), в похожем, но немного более громоздком (с большим числом параметров и оценивающих последовательностей). Опуская в правой части $V(x_*,x_{N+1})$, далее оптимально подбирают параметры метода $\{\lambda_k\}$, получают оценку скорости сходимости метода по функции в среднем. Если считать, что $\|x_*-x_k\|=\mathrm{O}(R)$, то отсюда также получают оценки скорости сходимости и с вероятностями больших уклонений (используется обобщение неравенства Азума–Хефдинга для последовательности мартингал-разностей [24, 79]). В детерминированном случае соотношение $\|x_*-x_k\|=\mathrm{O}(R)$ имеет место (всегда в евклидовом случае, и в зависимости от метода в общем случае) ввиду сходимости метода и того, что параметры оптимально подбираются так, что слагаемые $V(x_*,x_0)$ и $\tilde{\Delta}_N$ одного порядка (отличаются обычно не более чем в 10 раз):

$$\frac{1}{2}\|x_k-x_*\|^2\le V(x_*,x_k)\le V(x_*,x_0)+\tilde{\Delta}_{k-1}\le V(x_*,x_0)+\tilde{\Delta}_N .$$

В случае стохастического оракула, к сожалению, такие рассуждения уже не проходят. Можно, однако, из таких соображений оценить $E_{x_k}\left[V(x_*,x_k)\right]$. Дальше угадывается хвост распределения случайной величины $\|x_*-x_k\|$ исходя из выписанного выше соотношения, которое стоит понимать как равенство, т.е. хвост распределения ищется как неподвижная точка (а точнее ее оценка). Задавшись определенным доверительным уровнем $\sigma\ge 0$ можно оценить “эффективный” $R$: с вероятностью $\ge 1-\sigma$ имеют место неравенства $V(x_*,x_k)\le R$, $k=0,...,N$. В частности, для субгауссовских стохастических градиентов $R=\mathrm{O}\left(V(x_*,x_0)\ln^2(N/\sigma)\right)$, а для равномерно ограниченных – $R=\mathrm{O}\left(V(x_*,x_0)\ln(N/\sigma)\right)$. Детали можно посмотреть в доказательстве теоремы 4 работы [60] (см. также [43]). Примечательно, что все эти рассуждения в случае не ограниченного множества $Q$ не требуют равномерной ограниченности констант Липшица (функции, градиента) на всем $Q$ [92]. Похожим образом можно получать оценки вероятностей больших уклонений в сильно выпуклом случае в онлайн контексте (см. конец этого пункта). К сожалению, не все методы обладают такими же свойствами. Например, PGM [21] (в случае детерминированного оракула и не евклидовой прокс-структуры [43]) гарантирует лишь, что $\|x_*-x_k\|=\mathrm{O}(R)$, $k=0,...,N$, если [21] (прокс-диаметр здесь не нужен):

$$R=\max\left\{\|x-x_*\|:\ x\in Q,\ f(x)\le f(x_0)\right\}.$$

Хотя PGM и является релаксационным методом ( $f(x_{k+1}) \le f(x_k)$), возможно, что $R = \infty$. Требование $R < \infty$ (коэрцитивности) не является сильно обременительным. Его можно обеспечивать за счет регуляризации задачи [97].

Полезно также иметь в виду, что за счет допускаемой неточности оракула, можно погрузить задачу с гельдеровым градиентом, т.е. удовлетворяющим неравенству $\|\nabla f(x) - \nabla f(y)\|_* \le L_\nu \|x - y\|^\nu$, при некотором $\nu \in [0,1]$ (в том числе и негладкую задачу с ограниченной нормой разности субградиентов при $\nu = 0$) в класс гладких задач с неточным оракулом, характеризующимся точностью $\delta$ и [78]

(10) $$L = L_\nu \left[ \frac{L_\nu (1-\nu)}{2\delta(1+\nu)} \right]^{\frac{1-\nu}{1+\nu}}.$$

Заметим, в этой связи, что если в предположении 1 считать

$$E_\xi\left[f(y,\xi)\right] - E_\xi\left[F(x,\xi)\right] - \left\langle E_\xi\left[G(x,\xi)\right], y - x \right\rangle \le \frac{L}{2}\|y - x\|^2 + M\|y - x\| + \delta,$$

то вместо $D$ в (7) стоит писать $M^2 + D$ [95, 96].

Таким образом, например, можно получить оценки (2), (4) из оценки (7). В частности, метод двойственных усреднений и зеркальный спуск (см. п. 2) можно получить из PGM в варианте [79] с неточным оракулом и $L = M^2/(2\delta)$. Но наряду с введенной нами искусственной неточностью оракула, можно допустить, что имеется также реальная неточность оракула. Несложно привести оценки (на базе формулы (7) и ((8), (9))) сочетающие наличие в задаче искусственной и реальной неточности [61].

Как мы предполагали выше, множество $Q$ должно быть достаточно простой структуры, чтобы на него можно было эффективно проектироваться. Однако в приложениях часто возникают задачи условной минимизации [2], в которых, например, есть ограничения вида $g(x) \le 0$, где $g(x)$ – выпуклые функции [15]. "Зашивать" эти ограничения в $Q$, как правило, не представляется возможным в виду вышесказанного требования о легкости проектирования. Тем не менее, на основе описанного выше можно строить (за дополнительную логарифмическую плату) двухуровневые методы (наверное, лучше говорить "методы уровней", чтобы не было путаницы с многоуровневой оптимизацией, см. пример 5) условной оптимизации [10, 96]. При этом на каждом шаге такого метода потребуется проектироваться на пересечение множества $Q$ с некоторым полиэдром, вообще говоря, зависящим от номера шага. Последнее обстоятельство в общем случае сужает класс задач, к которому применимы такие многоуровневые методы до класса задач, к которым применимы методы внутренней точки [10]. В частности, возникает довольно обременительное условие на размер пространства, в котором проходит оптимизация: $n \sim 10^4 - 10^5$. Все это не удивительно, поскольку имеются нижние оценки [4] (рассматриваются аффинные ограничения в виде равенств, аналогично могут быть рассмотрены и неравенства), показывающие, что в общем случае для нахождения такого $x \in \mathbb{R}^n$, что $\|Ax-b\|_2 \le \varepsilon$ потребуется не меньше, чем $\Omega\left(\sqrt{L_x} R_x / \varepsilon\right)$ операций типа умножения $Ax$ ($L_x = \sigma_{\max}(A) = \lambda_{\max}\left(A^T A\right)$ с – максимальное собственное значение матрицы $A^T A$, $R_x = \left\|x^*\right\|_2 = \left\|\left(A^T A\right)^{-1} A^T b\right\|_2$ ). Аналогичное можно сказать и про седловые задачи: для отыскания такой пары $(x, y)$, что (левая часть этого неравенства всегда неотрицательная)

$$\max_{\tilde{y} \in S_n(1)} \tilde{y}^T A x - \min_{\tilde{x} \in S_n(1)} y^T A \tilde{x} \le \varepsilon$$

потребуется не меньше, чем $\Omega(\Lambda/\varepsilon)$ ( $\Lambda$ – максимальный по модулю элемент матрицы $A$) операций типа умножения $Ax$ и $y^T A$. Заметим, что обе выписанные нижние оценки справедливы при условии, что число итераций (операций типа умножения $Ax$) $k \le n$. Как следствие, в общем случае даже для гладкой детерминированной сильно выпуклой постановки при наличии всего лишь аффинных ограничений $Ax = b$ нельзя надеяться на быстрое решение. Тем не менее, некоторые дополнительные предположения в ряде случаев позволяют ускорить решение таких задач (см., например, [55, 97]).

**Замечание 5 (см. также [92]).** Задача поиска такого $x^*$, что $Ax^* = b$ сводится к задаче выпуклой гладкой оптимизации

$$f(x) = \|Ax-b\|_2^2 \to \min_x .$$

Нижняя оценка для скорости решения такой задачи [4] (см. также формулу (7) с $\delta = D = 0$, $p = 1$) имеет вид: $f(x_k) \ge \Omega\left(L_x R_x^2 / k^2\right)$. Откуда следует, что только при $k \ge \Omega\left(\sqrt{L_x} R_x / \varepsilon\right)$ можно гарантировать выполнение неравенства $f(x_k) \le \varepsilon^2$, т.е. $\|Ax_k - b\|_2 \le \varepsilon$.

Заметим, что эта нижняя оценка для специальных матриц может быть улучшена. Причем речь идет не о недавних результатах D. Spielman'а [98] (премия Неванлины 2010 года), а о более простой ситуации. Вернемся к задаче поиска вектора PageRank (примеры 1, 2), которую мы перепишем как

$$Ax = \begin{pmatrix} \left(P^T - I\right) \\ 1 \ldots\ldots 1 \end{pmatrix} x = \begin{pmatrix} 0 \\ 1 \end{pmatrix} = b \ , \ I \text{ – единичная матрица.}$$

По теореме Фробениуса–Перрона [99] решение такой системы с неразложимой матрицей $P$ единственно и положительно $x > 0$. Сведем решение этой системы уравнений к вырожденной задаче выпуклой оптимизации

$$\frac{1}{2}\|x\|_2^2 \to \min_{Ax=b} .$$

Построим двойственную к ней задачу [9]

$$\min_{Ax=b} \frac{1}{2}\|x\|_2^2 = \min_x \max_\lambda \left\{\frac{1}{2}\|x\|_2^2 + \langle b - Ax, \lambda\rangle\right\} = \max_\lambda \min_x \left\{\frac{1}{2}\|x\|_2^2 + \langle b - Ax, \lambda\rangle\right\} = \max_\lambda \left\{\langle b, \lambda\rangle - \frac{1}{2}\left\|A^T \lambda\right\|_2^2\right\}.$$

Поскольку система $Ax = b$ совместна, то по теореме Фредгольма не существует такого $\lambda$, что $A^T\lambda = 0$ и $\langle b, \lambda \rangle > 0$, следовательно, двойственная задача имеет конечное решение (т.е. существует ограниченное решение двойственной задачи $\lambda^*$). Зная решение $\lambda^*$ двойственной задачи

$$\langle b, \lambda \rangle - \frac{1}{2}\left\| A^T \lambda \right\|_2^2 \to \max_{\lambda}$$

можно восстановить решение прямой задачи (из условия оптимальности по $x$): $x(\lambda) = A^T \lambda$. Однако важно здесь то, что FGM [10] для этой двойственной задачи дает возможность попутно получать следующую оценку на норму этого градиента [92]:

$$\left\| Ax_k - b \right\|_2 = \mathrm{O}\left( \frac{L_\lambda R_\lambda}{k^2} \right),$$

где $x_k$ есть известная выпуклая комбинация

$$\left\{ x(\lambda_i) \right\}_{i=1}^k, \;\; L_\lambda = \sigma_{\max}\left(A^T\right) = \sigma_{\max}(A), \;\; R_\lambda = \left\| \lambda^* \right\|_2,$$

где можно считать, что $\lambda^*$ – решение двойственной задачи с наименьшей евклидовой нормой. Кажется, что это противоречит нижней оценке $\left\| Ax_k - b \right\|_2 \geq \Omega\left( \sqrt{L_x} R_x / k \right)$. Однако, важно напомнить [4], что эта нижняя оценка установлена для всех $k \leq n$ ($n$ – размерность вектора $x$), и она будет улучшена, в результате описанной процедуры только если дополнительно предположить, что матрица $A$ удовлетворяет следующему условию $L_\lambda R_\lambda \ll n\sqrt{L_x} R_x$, что сужает класс, на котором была получена нижняя оценка $\Omega\left( \sqrt{L_x} R_x / k \right)$. В типичных ситуациях можно ожидать, что $R_\lambda \gg R_x$ ($R_x \leq \sqrt{2}$). Это обстоятельство мешает выполнению требуемого условия.

**Пример 4.** Если имеется дополнительная информация о структуре седловой задачи, то можно её использовать для ускорения [14, 100]. Более того, многие современные постановки задач (негладкой) выпуклой оптимизации (в частности, связанные с compressed sensing и $l_1$-оптимизацией) в пространствах огромных размеров специально стараются представить седловым образом с целью получения эффективного решения (см. работы А.С. Немировского, А.Б. Юдицкого [31, 101, 102]). Далее будет разобран один простой пример (немного обобщающий результаты [78, 79], см. также [92]), демонстрирующий возможности градиентных методов с неточным оракулом в седловом контексте. Рассматривается седловая задача ($x \in \mathbb{R}^n, y \in \mathbb{R}^m$)

$$f(x) = \max_{\|y\|_2 \leq R_y} \left\{ G(y) + \langle By, x \rangle \right\} \to \min_{\|x\|_2 \leq R_x},$$

где функция $G(y)$ – сильно вогнутая с константой $\kappa$ относительно 2-нормы и константой Липшица градиента $L_G$ (также в 2-норме). Тогда функция $f(x)$ будет гладкой, с константой Липшица градиента в 2-норме $L_f = \sigma_{\max}(B)/\kappa$. Казалось бы, что мы можем решить задачу минимизации функции $f(x)$ за $\mathrm{O}\left( \sqrt{\sigma_{\max}(B) R_x^2 / (\kappa\varepsilon)} \right)$ итераций, где $\varepsilon$ - желаемая точность по функции. Но это возможно, только если мы можем абсолютно точно находить $\nabla f(x) = By^*(x)$, где $y^*(x)$ – решение вспомогательной задачи максимизации по $y$ при заданном $x$. В действительности, мы можем решать эту задачу (при различных $x$) лишь приближенно. Если мы решаем вспомогательную задачу быстрым градиентным методом [14] с точностью $\delta/2$ (на это потребуется $O\left(\sqrt{L_G/\mu}\ln\left(L_G R_y^2/\delta\right)\right)$ итераций), то пара $\left(G\left(y_{\delta/2}(x)\right) + \left\langle By_{\delta/2}(x), x\right\rangle, By_{\delta/2}(x)\right)$, где $y_{\delta/2}(x)$ – $\delta/2$-решение вспомогательной задачи, будет $\left(\delta, 2L_f, 0\right)$-оракулом [78, 79]. Выбирая $\delta = O\left(\varepsilon\sqrt{\varepsilon/\left(L_f R_x^2\right)}\right)$ (см. формулу (9) при $p = 1$), получим после

$$O\left(\sqrt{\frac{L_G \sigma_{\max}(B) R_x^2}{\kappa^2 \varepsilon}} \ln\left(\frac{L_f L_G R_x^2 R_y^2}{\varepsilon}\right)\right)$$

итераций (на итерациях производится умножение матрицы $B$ на вектор/строчку и вычисление градиента $G(y)$) $\varepsilon$-решение задачи минимизации $f(x)$. Отметим, что если не использовать сильную вогнутость функции $G(y)$, то для получения пары $(x_N, y_N)$, удовлетворяющей неравенству

$$\max_{\|y\|_2 \le R_y}\left\{G(y) + \langle By, x_N\rangle\right\} - \min_{\|x\|_2 \le R_x}\left\{G(y_N) + \langle By_N, x\rangle\right\} \le \varepsilon,$$

потребуется $\Omega\left(\max\left\{L_G R_y^2, \sigma_{\max}(B) R_x R_y\right\}/\varepsilon\right)$ итераций (см., например, [4, 9, 11]).

Интересно отдельно разобрать ситуацию, когда вместо множества $\|y\|_2 \le R_y$ фигурирует симплекс $S_m(R_y)$, $G(y) = -\sum_{k=1}^{m} y_k \ln\left(y_k/R_y\right)$ – сильно вогнутая в 1-норме с константой $\kappa = 1$ функция и $R_x = \infty$ (энтропийно-линейное программирование [97]). В этом случае мы не можем обеспечить даже равномерной ограниченности градиента функции $G(y)$. Тем не менее, также можно рассчитывать [97] на зависимость $O\left(\varepsilon^{-1/2}\ln\left(\varepsilon^{-1}\right)\right)$ числа итераций от точности $\varepsilon$ для критерия:

$$\max_{y \in S_m(R_y)} \min_{x}\left\{G(y) + \langle By, x\rangle\right\} - \min_{x}\left\{G(y_N) + \langle By_N, x\rangle\right\} \le \varepsilon .$$

При этом вместо энтропии в качестве функции $G(y)$ можно брать любую сильно вогнутую в 1-норме функцию, для которой решение задачи максимизации (вычисление $f(x)$ с точностью $\varepsilon$) может быть осуществлено за $O\left(\ln\left(\varepsilon^{-1}\right)\right)$. В примере с энтропией, для $f(x)$ есть просто явная формула. Точнее, важно то, что есть явная формула[27] для оптимального решения $y^*(x)$.[28] К сожалению, имеется проблема вхождения в оценку необходимого числа итераций неизвестного размера решения $x_*$ задачи минимизации $f(x)$. Эта проблема решаема [97]. В частности, в случае, когда $G(y)$ имеет ограниченную вариацию на множестве $S_m(R_y)$ (для энтропии эта вариация равна $R_y \ln m$), можно предложить метод, с оценкой числа итераций $O(\varepsilon^{-1}\ln(\varepsilon^{-1}))$. В эту оценку уже никак не входит неизвестный размер решения $x_*$, который может оказаться большим [97]. Далее мы еще вернемся к вопросу о том, как действовать, в случае, когда тот или иной параметр задачи (в данном случае размер решения) априорно не известен.

Отметим, что сильной вогнутости можно добиться и искусственно [77], глава 3 [79], [92, 97]. Подход отмеченных работ приводит к оптимальным для такого класса задач оценкам (с точностью до логарифмического фактора[29] $\ln(\varepsilon^{-1})$), и позволяет, на самом деле, контролировать точность решения одновременно по $x$ и по $y$ без использования прямо-двойственности в классическом варианте (см. ниже), что может быть полезным в определенных ситуациях [55]. Здесь под оптимальными методами мы имеем в виду методы с проксимальным оракулом. Однако в ряде задач оптимизации огромных размеров оказывается эффективнее использовать линейный минимизационный оракул [103], пришедший из метода Франк–Вульфа (см., например, п. 3.3 [11]). Грубо говоря, суть подхода в том, что сначала вычисляется не $f(x)$ согласно модели, описанной в примере 4, а в седловом представлении задачи меняется порядок взятия максимума и минимума, и вычисляется с помощью линейного минимизационного оракула сначала минимум по $x$. Причем это не обязательно делать точно (см. п. 5 § 1 главы 5 [2]). Получающаяся задача максимизации по $y$ уже не будет гладкой, поэтому с учетом сильной вогнутости $G(y)$ здесь можно рассчитывать только на зависимость числа итераций от желаемой точности $O\left(\varepsilon^{-1}\right)$. Получается вроде как хуже, чем раньше. Но тут надо учитывать, как входят размерности $n$ и $m$, которые могут быть огромными в приложениях, см. п. 3.3 [11], [59, 60, 61, 103]. Удивительным образом, в сложность внутренней задачи при таком подходе (минимизации по $x$) при определенной структуре (как правило, связанной с ограничениями симплексного типа и матрицей $B$, имеющей комбинаторую [103] или сетевую природу [60]) может не входить размерность вектора $x$ (т.е. $n$), что позволяет решать задачи колоссальных размеров по $n$.

[27] Сложность формулы оценивается числом ненулевых элементов в матрице $B$. При этом считаем, что градиент $G(y)$ рассчитывается быстрее, чем занимает умножение матрицы $B$ на столбец / строку.

[28] Отметим, что если $G(y)$ – сепарабельная вогнутая функция (но не обязательно сильно вогнутая) и вместо ограничения $\|y\|_2 \le R_y$ задано сепарабельное ограничение (например, $\|y\|_\infty \le R_y$), то $\varepsilon$-приближенный поиск $y^*(x)$ можно осуществить за $O(\ln(\varepsilon^{-1}))$ умножений матрицы $B$ на столбец, решая соответствующие одномерные задачи. Не много более громоздкие рассуждения [55] позволяют и при наличии ограничения $\|y\|_2 \le R_y$ осуществить $\varepsilon$-приближенный поиск $y^*(x)$ также за $O(\ln(\varepsilon^{-1}))$ умножений матрицы $B$ на столбец. К сожалению, отсутствие сильной вогнутости не позволяет использовать в том же виде концепцию $(\delta, L, \mu)$-оракула для внешней задачи, однако можно при этом использовать концепцию $\delta$-субградиента для внешней задачи [2]. Это приводит лишь к оценкам $O(\varepsilon^{-2})$, которые уже не будут оптимальными (улучшаемы до $O(\varepsilon^{-1})$).

[29] Ниже мы обсудим, как можно избавиться от этого логарифмического фактора для задач с явной формулой для $y^*(x)$, например, для задач энтропийно-линейного программирования [55, 97].

Пример 4 был приведен, прежде всего, потому, что он поясняет одно интересное и достаточно современное направление в численных методах выпуклой оптимизации (см., например, [61, 62, 92, 96, 101, 102]). Грубо говоря, это направление можно охарактеризовать, как попытку ввести оптимальную “алгебру” над алгоритмами выпуклой оптимизации. А именно, если требуется оптимизировать функционал (искать седловую точку), который обладает разными свойствами (гладкости, сильной выпуклости, быстроты вычислимости частных производных и т.п.) по разным группам переменных (такие задачи часто в последнее время возникают в разных приложениях, в частности, в транспортных и экономических [58, 60, 61, 62, 104–110]) и(или) сам представляет собой некоторую суперпозицию других функционалов (с разными свойствами; наиболее популярен случай суммы двух функционалов [9, 14, 43, 55, 59, 63, 96, 101, 102]), то хотелось бы получить такую декомпозицию исходной задачи, чтобы правильное сочетание (правильное чередование с правильными частотами) оптимальных методов для получившихся отдельных подзадач позволило бы получить оптимальный метод для исходной задачи. В ряде интересных случаев такое оказывается возможным (с оговоркой, что оптимальность понимается с точностью до логарифмического фактора). По-видимому, новым в этом направлении является наблюдение, отмеченное в примере 4 (см. также [61, 62, 92]), что при определенных условиях идея оптимального сочетания различных методов для решения одной сложной по структуре задачи оптимизации, может быть реализована на основе концепции неточного оракула.

Другой способ борьбы с дополнительными ограничениями типа равенств или неравенств в задачах выпуклой оптимизации базируется на прямо-двойственной структуре [27] всех обсуждаемых методов (поскольку они строят модель функции [14]). Это означает, что ограничения вносятся во вспомогательную задачу оптимизации, возникающую на каждом шаге метода и отвечающую за проектирование. В результате на каждом шаге решается более сложная задача. Тем не менее, если такие вспомогательные задачи можно эффективно решать (что в общем случае также наталкивается на сложности, описанные ранее) с помощью метода множителей Лагранжа (найдя и сами множители) или когда у исходной задачи есть модель (см. пример 4 и [18, 55, 59–62, 77, 92, 97, 110, 111]), то тогда описанные методы позволяют не только эффективно решать исходную задачу оптимизации с ограничениями, но и находить попутно (по явно выписываемым формулам) решение двойственной задачи.

Основная идея работы [27] состоит в том (здесь мы ограничимся рассмотрением детерминированного случая с точным оракулом, выдающим градиент; в стохастическом случае с неточным оракулом см., например, [55] и лемму 7.7 [79]), что метод генерирует в прямом пространстве на итерациях такую последовательность $\{x_k\}$,[30] что зазор двойственности (duality gap) $\Delta(\lambda, x; N)$ удовлетворяет условию

$$\Delta(\lambda, x; N) = \max_{u \in Q} \left\{ \frac{1}{S_N} \sum_{k=0}^{N} \lambda_k \left\langle \nabla f(x_k), x_k - u \right\rangle \right\} \le \varepsilon ,$$

где $S_N = \sum_{k=0}^{N} \lambda_k$ , $\lambda_k \ge 0$ , поэтому

$$f\left( \frac{1}{S_N} \sum_{k=0}^{N} \lambda_k x_k \right) - \min_{x \in Q} f(x) \le \varepsilon .$$

Это следует из выкладки

$$f\left( \frac{1}{S_N} \sum_{k=0}^{N} \lambda_k x_k \right) - f(u) \le \frac{1}{S_N} \sum_{k=0}^{N} \lambda_k \cdot \left( f(x_k) - f(u) \right) \le \frac{1}{S_N} \sum_{k=0}^{N} \lambda_k \left\langle \nabla f(x_k), x_k - u \right\rangle .$$

Аналогичную точность (для двойственной задачи) дает следующая аппроксимация решения двойственной задачи

$$y = \frac{1}{S_N} \sum_{k=0}^{N} \lambda_k y_k .$$

Это сразу следует из того, что зазор двойственности оценивает сверху разность между получившимися значениями целевой функции в прямой задаче и двойственной [27], которую мы будем называть двойственным зазором. Эта разность всегда неотрицательна, и на точных решениях прямой и двойственной задачи (и только на них) равна нулю. Заметим, что контроль онлайн-зазора двойственности

$$\Delta(\lambda, x; N) = \max_{u \in Q} \left\{ \frac{1}{S_N} \sum_{k=0}^{N} \lambda_k \left\langle \nabla f_k(x_k), x_k - u \right\rangle \right\}$$

позволяет в случае, когда удается выбрать $\lambda_k \equiv 1$, получать оценки регрета (псевдо регрета в стохастическом случае) в задачах онлайн оптимизации (см. [32, 45] и конец этого пункта). К сожалению, ограничение $\lambda_k \equiv 1$ существенно сужает класс методов. Скажем, для рассматриваемых в этом пункте быстрых градиентных методов $\lambda_k \sim k^p$. Кроме того, даже если в онлайн постановке допустить взвешивание с различными весами, все равно требуется, чтобы способ получения оценки на зазор двойственности допускал бы обобщение на онлайн-постановки. Быстрый градиентный метод, например, этого не допускает, что не сложно усмотреть из оценок работы [21].

[30] В двойственном пространстве при этом генерируется последовательность соответствующих множителей Лагранжа $\{y_k\}$ [27] или, в случае наличия модели у исходной прямой задачи (см. пример 4), последовательность $\{y_k\}$ генерируется по явным или расчетным формулам $\{y_k = y(x_k)\}$ согласно этой модели [18, 55, 59, 92, 110, 111]. Такой подход также позволяет убрать логарифмический фактор в задачах энтропийно-линейного программирования [55, 97] и аналогичных задачах, см. п. 5.2 [111] и [43, 55, 92, 112].

Описанная выше конструкция, основанная на оценке зазора двойственности, работает в случае ограниченного множества $Q$. В случае неограниченного $Q$ (это типичная ситуация, когда необходимо решать двойственную задачу, по решению которой требуется восстанавливать решение прямой задачи) можно искусственно компактифицировать $Q$ [55, 60, 61, 111]. Однако, в большинстве случаев такая компактификация не позволяет очевидным образом оценивать настоящий (не обрезанный) двойственный зазор в исходной задаче, что часто представляется важным ввиду наличия простых явных формул для этого настоящего зазора, и возможности использования контроля зазора двойственности у качестве критерия останова метода. Несмотря на отмеченную теоретическую проблему, на практике проблема оказывается решаемой [55, 59, 60].

Более общий способ оценки разности между получившимися значениями целевых функций в прямой задаче и двойственной базируется на контроле сертификата точности [111, 113] (accuracy certificate), в который наряду с градиентами функционала входят градиенты нарушенных ограничений или в общем случае вектора нормалей к гиперплоскостям, отделяющим $x_k$ от множества $Q$ (в ряде постановок "градиенты" стоит заменить на "субградиенты"). Вектора двойственных множителей формируются из соответствующих (сертификату точности) взвешенных сумм векторов нормалей отделяющих гиперплоскостей [111, 113, 114]. Собственно, такая интерпретация двойственных множителей следует из способа обоснования принципа множителей Лагранжа на основе следствия теоремы Хана–Банаха (теоремы об отделимости) [115]. Причем в работе [113, 114] за счет слейтеровской релаксации ограничений (допущения возможности нарушения ограничений на $\varepsilon$ [9, 97]) получаются оценки скорости сходимости, не зависящие от размера двойственного решения, который может быть большим.

Во многих (транспортно-)экономических приложениях при поиске равновесных конфигураций требуется решать пару задач (прямую и двойственную), см., например, [58–62, 97, 104–110]. Причем интересны решения обеих задач (решения этих задач имеют содержательную интерпретацию и используются при принятии решений / управлении). Если у этой пары задач, на которую можно смотреть, как на одну седловую задачу, есть определенная структура (проявляющаяся, например, в сильной выпуклости функционала по части переменных, наличии эффективно вычислимого линейного минимизационного оракула и т.п.), то описанный выше формализм позволяет развить идею примера 4 таким образом, чтобы одновременно (без дополнительных затрат) получать решения обеих задач. Даже в случае огромного размера одной из этих задач, можно надеяться (при эффективном линейном минимизационном оракуле), что эта размерность не войдет в сложность поиска решения прямой и двойственной задачи [59, 60, 103, 110].

В действительности, выбранный в данной статье класс проекционных методов с построением модели функции далеко не единственный возможный способ строить прямо-двойственные методы. Скажем, уже упоминавшиеся методы условного градиента также являются прямо-двойственными [17, 18]. Еще более удивительным может показаться, что прямо-двойственная интерпретация есть, например, у метода эллипсоидов [111]. Более того, в ряде ситуаций мы можем за линейное время (с геометрической скоростью сходимости) находить одновременно решение прямой и двойственной задачи. Причем речь идет не только о конструкциях типа [74], базирующихся на принципе (см. также замечание 6): сопряженная функция к выпуклой функции с липшицевым градиентом – сильно выпуклая, и обратно, сопряженная к сильно выпуклой функции – выпуклая функция с липшицевым градиентом; но и о более общем контексте [43, 55, 74, 92, 111].

Возвращаясь к сказанному выше в связи с оценками (7) – (10) интересно заметить, что если множество $Q \subset \mathbb{R}^n$ есть шар $B_q^n(R)$ радиуса $R$ в $q$-норме,[31] то нижние оценки (для случая $D=\delta=0$) на точность (по функции), которую можно получить после $N \le n$ итераций, имеют вид [116] (считаем, что[32] $\|\nabla f(x)-\nabla f(y)\|_{q'} \le L_\nu \|x-y\|_q^\nu$, $1/q+1/q'=1$, $\nu \in (0,1]$):

---

[31] Если о выпуклом замкнутом множестве $Q$ известно только то, что оно содержит $B_q^n(R)$, то все сказанное далее также остается в силе.

[32] В случае достаточной гладкости функции $f(x)$ можно выписать следующее представление для константы Липшица градиента (верхний индекс $q$ соответствует выбору нормы в прямом пространстве)

$$L^q = \max_{x\in Q, \|h\|_q \le 1} \langle h, \nabla^2 f(x) h \rangle .$$

В частности, $L^1 \le L^2 \le nL^1$, $L^2 \le L^\infty \le nL^2$. Эти формулы вместе со сказанным ранее относительно того, как может меняться $M$ – обычная константа Липшица $f(x)$, при изменении нормы в прямом пространстве, поясняют почему в “устойчивые сочетания” эти константы входят таким образом: $M^2R^2/\varepsilon^2$, $LR^2/\varepsilon$. Если ввести “физические размерности”, скажем, считать, что $f(x)$ это рубли (*руб*), а $x$ это килограммы (*кг*), то $\varepsilon$ [*руб*], $R$ [*кг*], $M$ [*руб*/*кг*], $L$ [*руб*/*кг*$^2$]. Поскольку число итераций $N$ должно быть безразмерной величиной, то возникновение агрегатов $M^2R^2/\varepsilon^2$, $LR^2/\varepsilon$ вполне закономерно. Аналогичные рассуждения можно провести и для оценок в сильно выпуклом случае. Все это приводит к довольно интересным следствиям [92]. Например, что шаг метода в негладком случае $h \sim \varepsilon/M^2$, в гладком случае $h$ определяется из соотношения вида ($W(\ )$, $\tilde{W}(\ )$ – какие-то функции)

$$W\left(h\frac{M^2}{\varepsilon}, hL\right) = 1.$$

В стохастическом случае (вместо градиента получаем стохастический градиент с дисперсией $\sigma^2$) из

$$\Omega\left(\frac{1}{\min\{q,\ln n\}^{\nu}}\frac{L_{\nu}R^{1+\nu}}{N^{\nu+(\nu+1)/q}}\right)\ (2\le q\le\infty),$$

$$\Omega\left(\frac{1}{\ln^{\nu}(N+1)}\frac{L_{\nu}R^{1+\nu}}{N^{\nu+(\nu+1)/2}}\right)\ (1\le q<2).$$

Приведенный результат хорошо соответствует тому, что написано в замечании 2 (см. также [117]).

Для $q=\infty$ и $\nu=1$ (гладкий случай) приведенная оценка с точностью до логарифмического фактора будет иметь вид $\Omega\left(LR^2/N\right)$. Эта оценка достигается, например, на методе условного градиента Франк–Вульфа [17, 118].[33] Исходя из только что написанного и тезиса о неулучшаемости оценок (7) ((8), (9)) (при $D=\delta=\mu=0$, $p=1$) может возникнуть ощущение противоречия. Это ощущение дополнительно усиливается примером 2 из п. 2. Действительно, исходя из этого примера, может сложиться ощущение, что проблема выбора прокс-структуры в задаче не очень актуальна, поскольку можно исходить просто из самой нормы. И это действительно так, если мы ограничиваемся не ускоренными градиентными методами (PGM, метод Франк–Вульфа), которые сходятся как $O\left(L^qR_q^2/N\right)$ (здесь $R_q=R$ – диаметр множества $Q$, посчитанный в $q$-норме, в нашем случае $q=\infty$). Если же мы хотим ускориться, и достичь оптимальной оценки $O\left(LR^2/N^2\right)$, то уже необходимо существенно использовать прокс-функцию $d(x)\ge 0$ со свойством сильной выпуклости относительно выбранной нормы и с константой сильной выпуклости $\alpha\ge 1$ [14]. Скажем (в связи с примером 2), квадрат 1-нормы – не есть сильно выпуклая функция относительно 1-нормы, т.е. $d(x)=\|x\|_1^2$ – не может быть прокс-функцией при выборе 1-нормы для симплекса. В классе удовлетворяющих условию 1-сильной выпуклости прокс-функций (относительно выбранной нормы) подбирается такая, которая минимизирует $R^2:=\max_{x\in Q}d(x)$. Именно это $R^2$ входит в оценку FGM $O\left(LR^2/N^2\right)$. И как уже отмечалось (см. п. 2) для $Q=B_\infty^n(R)$ имеет место следующая оценка на прокс-диаметр $R^2:=\max_{x\in Q}d(x)=R^2\Omega(n)$. Отсюда, с учетом того, что $N\le n$, получаем, что

---

$$\tilde{W}\left(h\frac{M^2}{\varepsilon},hL,h\frac{\sigma}{R}\right)=1.$$

[33] Отметим также, что этот метод допускает обобщение на случай неточного оракула, и неулучшаемость оценок может быть проинтерпретирована с точки зрения сохранения свойства разреженности решения [17]. Это неудивительно, поскольку аналогичный метод (с линейным минимизационным оракулом, см. пример 4) с аналогичными оценками скорости можно получить (см. п. 5.5.1 [9], [57]) из композитного варианта FGM в концепции неточного оракула (суть метода в том, что в композитном варианте FGM на каждой итерации решается задача, в которой коэффициент при прокс-слагаемым равен нулю, т.е. оно просто отсутствует).

оценка $O\left(LR^2/N\right)$ и оценка $O\left(LR^2n/N^2\right)$, приводят, в общем-то, к одному результату, но в случае использования FGM требуется дополнительно искать оптимальную прокс-структуру. Только в таком случае будет совпадение результатов. Более того, также как и в замечании 2, здесь хорошо видно, что при $q \geq 2$ можно ограничиться рассмотрением только евклидовой прокс-структуры для FGM и евклидовой нормы для метода Франк–Вульфа. В частности, для $Q = B_\infty^n(R)$ действуя так, мы получим для FGM оценку $O\left(L^2R^2n/N^2\right)$ вместо ранее полученной оценки $O\left(L^\infty R^2n/N^2\right)$, соответствующей при $N \leq n$ неулучшаемой оценке $O\left(L^\infty R^2/N\right)$ (здесь мы проставили верхние индексы у $L$, поскольку они различаются). Дальше можно написать все тоже по поводу неулучшаемости оценок, что и в конце замечания 2.

Отметим также, что параметры $R$ и $\mu$ могут быть не известны априорно или процедуры их оценивания приводят к слишком (соответственно) завышенным и заниженным результатам. Это может быть проблемой, поскольку в ряде случаев знание этих и других параметров требуется методу для расчета величин шагов и условий остановки. Из этой ситуации можно выйти за логарифмическое (по этим параметрам) число рестартов метода. Стартуя, скажем, с $R = 1$ и делая число шагов, вычисленное из оценки скорости сходимости при выбранном $R$, мы проверяем выполняется ли для вектора, выдаваемого алгоритмом, условие $\varepsilon$-близости по функции (при условии, что мы можем сделать такую проверку). Если условие $\varepsilon$-близости не выполняется, то полагаем $R := 2R$ и т.д. Все эти перезапуски увеличат общее число обращений к оракулу лишь в $O(1)$ раз [14, 31, 38].[34] Аналогичное можно сказать про[35] $L$, $M$, $\mu$ и $D$. Однако если убрать стохастичность, тогда $L$, $M$ можно не только эффективнее адаптивно подбирать (аналогично правилу Армихо [16] в независимости от того можем ли мы сделать проверку условия $\varepsilon$-близости значения функции в текущей точке к оптимальному) по ходу самих итераций (увеличив в среднем число обращений к оракулу не более чем в 4 раза), но и в некотором смысле оптимально самонастраиваться (используя формулу (10)) на гладкость функционала на текущем участке пребывания метода [119]. Это означает, что в детерминированном случае без учета сильной выпуклости функционала существует универсальный метод, работающий по оценкам (8) с $L$, рассчитанной по формуле (10) (в которой $\delta$ берется из (9)), и оптимальным в смысле скорости сходимости выбором параметра $\nu \in [0,1]$. Причем выбор $\nu$ осуществляется не нами заранее, исходя из знания всех констант и минимизации выписанных оценок, а самим алгоритмом (здесь выбрано $p=1$):

$$N(\varepsilon) = \inf_{\nu \in [0,1]} \left( \frac{2^{\frac{3+5\nu}{2}} L_\nu R^{1+\nu}}{\varepsilon} \right)^{\frac{2}{1+3\nu}}.$$

---

[34] В свою очередь можно поиграть и на этом $O(1)$, стараясь его минимизировать. Для этого шаг, который мы для простоты положили равным 2, подбирают оптимально исходя из того, с каким показателем степени входит неизвестный (прогоняемый) параметр в оценку числа итераций [56, 97].

[35] Впрочем, в детерминированных постановках мы можем явно наблюдать за последовательностью выдаваемых оракулом субградиентов и отслеживать условие на норму субградиентов. Как только наше предположение нарушилось (при этом мы не успели сделать предписанное текущему $M$ число шагов), мы увеличиваем $M$ в два раза и перезапускаем весь процесс с новым значением $M$. Число таких перезапусков будет не более чем логарифмическим от истинного значения $M$. Все эти рассуждения с небольшими оговорками (типа равномерной п.н. ограниченности стохастического субградиента) переносятся и на стохастические постановки, в которых наблюдается стохастический субградиент. Для определенного класса задач, в которые неизвестные параметры входят только в критерий останова, но не в сам метод (к таким задачам, например, относится задача поиска равновесия в модели Бэкмнана методом Франк–Вульфа и неизвестной константе $L$) можно обходиться и без перезапусков [60]. Заметим также, что у ряда популярных методов (например, метода зеркального спуска) есть варианты, в которые входит не оценка супремума нормы субградиента (или градиента), а норма субградиента на текущей итерации, которая известна [9, 27, 60, 110].

Это соответствует (с точностью до логарифмического фактора) нижним оценкам [111], выписанным выше для случая $q \in [1,2]$. Отметим, что здесь при определении $R$ используется соответствующая прокс-функция, см. замечание 2. К сожалению, пока не очень понятно можно ли что-то похожее сделать с параметром $\mu$ и с введенным нами в начале этого пункта параметром метода $p \in [0,1]$. Обзор других работ на тему самонастройки алгоритмов в гладком детерминированном случае имеется в [20, 120, 121], а в стохастическом случае в [122].

Приведенную оценку можно обобщить, если дополнительно известно, что функция $f(x)$ – $\mu$-сильно выпукла. Можно дополнительно к искусственно введенной игре на неточности оракула допустить, что имеет место и настоящая неточность. В этом случае также можно выписать соответствующие оценки [61]. Не играя на выборе $\nu \in [0,1]$, можно распространить все, что описано выше в этом абзаце на стохастические постановки. Аналогичное можно сделать для стохастических безградиентных и покомпонентных методов с неточным оракулом (см. п. 4 и [43]). Соответствующие обобщающие формулы собраны в работе [120–122], мы не будем их здесь приводить. Такие обобщения востребованы, например, в связи с приложениями к поиску равновесий в многостадийных моделях равновесного распределения транспортных потоков [58, 59, 61, 62, 107–109]. В основе этих приложений лежит конструкция, изложенная в примере 4, с универсальным методом [119] вместо FGM для решения внешней задачи.

Выше мы сделали обременительное предположение о возможности выполнять проверку условия $\varepsilon$-близости по функции. Такое заведомо возможно только при известном значении функционала в точке оптимума. Как правило, такой информации у нас априорно нет. Один из способов выхода из этой ситуации для задач стохастической оптимизации описан в п. 7.7 работы [79]. Другой способ – контролировать зазор двойственности (со стохастическими градиентами). Для применимости этого способа требуется, чтобы числовая последовательность $\left\{\lambda_k/S_N\right\}_{k=0}^{N}$ не зависела от неизвестных параметров. Во многих задачах, приходящих из транспортных и экономических приложений, нужно одновременно находить решения прямой и двойственной задачи, которые можно явно выписать. В таких случаях имеется эффективный способ проверки условия $\varepsilon$-близости по функции. Нужно проверить условие $\varepsilon$-малости разницы между полученным (приближенным) значением функционалов прямой и двойственной задачи (т.е. двойственного зазора) [59–61, 109].

Отметим также, что в детерминированном $\mu$-сильно выпуклом случае, когда в точке минимума $x_*$ выполняется условие[36] $\nabla f\left(x_*\right)=0$, критерий $\varepsilon$-близости по функции может быть переписан в терминах малости рассчитываемого на итерациях градиента:

$$f\left(x_k\right)-f_*\le\frac{1}{2\mu}\left\|\nabla f\left(x_k\right)\right\|_*^2\le\varepsilon .$$

В постановках с сильно выпуклой/вогнутой двойственной задачей (этого можно добиться искусственно, вводя регуляризацию в двойственную задачу, см. главу 3 [79] и [92, 97]) также можно оценивать точность решения прямой задачи по точности решения двойственной задачи, применяя к двойственной задаче неравенство [37]

$$\frac{1}{2L}\left\|\nabla f\left(x_k\right)\right\|_*^2\le f\left(x_k\right)-f_* .$$

В частности, это обстоятельство используется в критерии остановки двойственного метода из [97], см. также пример 4 и [55, 92].

Все сказанное выше, по-видимому, переносится в полной мере на задачи композитной оптимизации [9, 14, 78, 124] и некоторые их обобщения, например, [96, 101, 102].

**Замечание 6.** Композитные задачи имеют вид: $f\left(x\right)+\lambda h\left(x\right)\to\min\limits_{x\in Q}$, где $\lambda>0$, $h\left(x\right)$ – выпуклая функция простой структуры, скажем $h\left(x\right)=\left\|x\right\|_1$. Хочется, чтобы сложность решения этой задачи всецело определялась только гладкостью выпуклого функционала $f\left(x\right)$, а сильная выпуклость – обоими слагаемыми. Если не лианеризовывать функцию $h\left(x\right)$ при подсчете на каждой итерации градиентного отображения [14], а просто оставлять это слагаемое как есть, то, конечно, сложность решения вспомогательной задачи на каждой итерации увеличится (впрочем, в виду простой структуры функции $h(x)$, ожидается, что не намного), зато в оценку необходимого числа итераций уже не будут входить никакие константы, характеризующие гладкость $h(x)$, только константы, характеризующие сильную выпуклость (если имеется). На такие задачи также можно смотреть следующим образом (принцип множителей Лагранжа): $f(x) \to \min_{x \in Q,\, h(x) \le C(\lambda)}$. Поскольку функция $h(x)$ простой структуры, то проектироваться на множество Лебега этой функции несложно. Отсюда можно усмотреть независимость числа итераций от $h(x)$. Другой способ “борьбы” с композитным членом $h(x)$ (А.С. Немировский) заключается в переписывании задачи в “раздутом” (на одно измерение) пространстве: $f(x) + y \to \min_{x \in Q,\, h(x) \le y}$. Норма в раздутом пространстве задается как $\|(x,y)\| = \|x\| + \alpha |y|$. Функционал имеет такой вид, что в независимости от гладкости $f(x)$ в оценки супремума нормы субградиента / константы Липшица градиента не будет входить что-либо, связанное с $y$. В гладком случае все ясно сразу из определения, а в случае негладкой $f(x)$ это связано с тем, что, в действительности, в оценку необходимого числа итераций входит не супремум нормы субградиента, а супремум нормы разностей субградиентов [25] (см. также начало п. 2). За счет возможности выбирать сколь угодно маленьким $\alpha > 0$, можно считать независящим от $y$ и прокс-расстояние (от точки старта до решения раздутой задачи), входящее в оценку необходимого числа итераций. Таким образом, можно сделать оценку числа итераций независящей от $y$ и $h(x)$. В связи с написанным выше полезно заметить, что гладкость (липшицевость градиента) и сильная выпуклость функционала являются взаимодвойственными друг к другу для задач безусловной оптимизации (константа Липшица градиента переходит в константу сильной выпуклости и наоборот, отсюда, кстати сказать, можно усмотреть, что оценка скорости сходимости для таких задач должны зависеть от отношения $L/\mu$, другой способ понять это – соображения “физической” размерности), что активно используется в приложениях, см., например, [43, 55, 73, 74]. Однако для задач условной оптимизации остается только один переход: двойственная (сопряженная) задача к сильно выпуклой – гладкая (см., например, [77] и замечание 3), обратное не верно даже в случае сильно выпуклой функции $h(x)$. Собственно, мы уже сталкивались с “неравноправностью” гладкости и сильной выпуклости. При рассмотрении универсального метода, мы отмечали, что на гладкость можно настраиваться адаптивно, чего нельзя сказать про сильную выпуклость. Замечание 6 немного проясняет (с учетом Лагранжева формализма) соотношения между этими свойствами задачи. Впрочем, до окончательного понимания, к сожалению, сейчас еще довольно далеко. Не ясно даже принципиальны ли эти различия или их можно в перспективе устранить. По-видимому, принципиальны, но строгого обоснования мы здесь не имеем.

[36] От этого условия можно избавиться, используя в приведенных далее формулах вместо градиента градиентное отображение [10, 55].

[37] В замечании 5 (см. также [92]) был приведен пример, когда $\left\|\nabla f\left(x_k\right)\right\|_2=O\left(k^{-2}\right)$. Из выписанного неравенства мы можем гарантировать лишь $\left\|\nabla f\left(x_k\right)\right\|_2^2=O\left(k^{-2}\right)$. Ситуацию можно улучшить, если регуляризовать функционал (см. конец п. 2 и [92, 97]), сделав его сильно выпуклым, и применить FGM [10, 14, 92, 123] к регуляризованной задаче, тогда $\left\|\nabla f\left(x_k\right)\right\|_2=O\left(\left(\ln k\right)^2/k^2\right)$ (если использовать, например, FGM с оценкой числа итераций $O\left(\sqrt{L/\mu}\left\lceil\ln\left(\mu R^2/\varepsilon\right)\right\rceil\right)$, $\mu\sim\varepsilon/R^2$ ). В негладком случае ситуация проще, см. [27].

Также нам видится, что сказанное выше переносится на седловые задачи и монотонные вариационные неравенства [4, 14]. Причем, речь идет не о том, что было описано в примере 4, а о том, как скажется неточность оракула на оптимальные методы для седловых задач и монотонных вариационных неравенств [4, 14]. Ответ, более менее, известен: неточность оракула не будет накапливаться на оптимальных методах (в отличие от задач обычной выпуклой оптимизации). Отметим, что концепцию неточного оракула еще необходимо должным образом определить[38] – предположение 1 нуждается в

[38] Например, для монотонного вариационного неравенства: найти такой $x \in Q$, что для всех $y \in Q$ выполняется $\langle g(y), y - x \rangle \ge 0$, достаточные условия на $(\delta, L)$-оракул будут иметь вид: для любых $x, y \in Q$

$$\langle g(y) - g(x), y - x \rangle \ge -\delta, \; \|g(y) - g(x)\|_* \le L\|x - y\| + \delta .$$

Вероятно, в ряде ситуаций эти условия можно ослабить (доказательства в этом случае нам не известны)

$$-\delta \le \langle g(y) - g(x), y - x \rangle \le L\|y - x\|^2 + \delta .$$

корректировке для данного класса задач. Отсутствие накопления неточностей связано с тем, что для таких задач оценка (7) будет оптимальна (с некоторыми оговорками) при $p = 0$. Другие $p \in (0,1]$ рассматривать не стоит (такие оценки просто не достижимы). Впрочем, пока про это имеются лишь частичные результаты [96, 101, 102].

В оценку числа итераций для достижения заданной точности решения в описанных методах не входит явно размерность пространства $n$. Это наталкивает на мысль (подобные мысли, по-видимому, впервые были высказаны и реализованы для класса обычных градиентных методов в кандидатской диссертации Б.Т. Поляка [1], см. также [2, 15]) о возможности использовать эти методы, например, в гильбертовых пространствах [16, 42]. Оказывается это, действительно, можно делать при определенных условиях (см., например, [120] в контексте использованных в данной работе обозначений). В частности, концепция неточного оракула позволяет привнести сюда элемент новизны, существенно мотивированный практическими нуждами – принципиальной невозможностью (в типичных случаях нет явных формул) решать с абсолютной (очень хорошей) точностью вспомогательную задачу на каждом шаге градиентного спуска. Например, решение такой вспомогательной задачи для класса задач оптимального управления со свободным правым концом приводит к двум начальным задачам Коши для систем обыкновенных дифференциальных уравнений (важно, чтобы СОДУ для фазовых переменных и сопряженных решались, скажем, методом Эйлера, на одной и той же сетке), которые необходимо решить для вычисления градиента функционала [16]. Однако, в действительности, почти все практически интересные задачи (за редким исключением, к коим можно отнести класс Ляпуновских задач [115]) в бесконечномерных пространствах не являются выпуклыми, поэтому здесь имеет смысл говорить лишь о поиске локальных минимумов (локальной теории) [125]. Если ограничиться неускоренными методами (например, PGM), то можно показать, что при весьма общих условиях эти методы могут быть использованы в гильбертовом пространстве в концепции неточного оракула и для невыпуклых (но гладких) функционалов, причем с аналогичными оценками скорости сходимости (отличие от выпуклого случая будет в том, что метод сходится лишь к стационарной точке (локальному экстремуму), в бассейне притяжения которой окажется точка старта). Заметим, что задачи оптимального управления можно численно решать, построив соответствующую (аппроксимирующую) задачу оптимального управления с дискретным временем, что приводит к конечномерным задачам, для решения которых можно использовать конечномерный вариант PGM в невыпуклом случае (с точным оракулом). Этот путь, как правило, и предлагается в большинстве пособий (см., например, [16]). Однако при таком подходе мы должны уметь (по возможности точно) решать сложную задачу оценки качества аппроксимации исходной задачи оптимального управления ее дискретным по времени вариантом. Более теоретически обоснованный способ рассуждений, по сути, приводящий к необходимости решать все те же конечномерные задачи, заключается в рассмотрении исходной задачи оптимального управления и ее решения бесконечномерным вариантом PGM в невыпуклом случае (с неточным оракулом). Неточность оракула существенна. Поскольку на каждой итерации этого градиентного метода необходимо решать две задачи Коши для СОДУ, что в общем случае можно сделать лишь приближенно, но с лучшим контролем точности, чем при подходе с дискретизацией задачи оптимального управления. Отметим, что во многих “физических” приложениях схема Эйлера имеет хорошие теоретические свойства сходимости (устойчивости). Связано это с тем, что на оптимальном режиме, как правило, наблюдается некоторая стабилизация поведения системы управления, что приводит к устойчивости якобиана прямой и обратной системы дифференциальных уравнений.

Покажем в заключение этого пункта, как приведенные результаты переносятся на задачи стохастической онлайн оптимизации. Для этого напомним вкратце, в чем состоит постановка задачи (см., например, [9, 26, 32, 45, 126–131]). Требуется подобрать последовательность[39] $\{x^k\} \in Q$ так, чтобы минимизировать псевдо регрет:

$$\frac{1}{N}\sum_{k=1}^{N} E_{\xi^k}\left[f_k\left(x^k,\xi^k\right)\right] - \min_{x\in Q}\frac{1}{N}\sum_{k=1}^{N} E_{\xi^k}\left[f_k\left(x,\xi^k\right)\right]$$

на основе доступной информации $\left\{\nabla\tilde{f}_1\left(x^1,\xi^1\right);\ldots;\nabla\tilde{f}_{k-1}\left(x^{k-1},\xi^{k-1}\right)\right\}$ при расчете $x^k$, где[40]

$$\left\|\nabla\tilde{f}_k\left(x^k,\xi^k\right) - \nabla f_k\left(x^k,\xi^k\right)\right\|_* \le \delta,$$

$$E_{\xi_k}\left[\nabla f_k\left(x,\xi^k\right)\right] = \nabla f_k(x).$$

Здесь с.в. $\{\xi^k\}$ могут считаться независимыми одинаково распределенными. Онлайновость постановки задачи допускает, что на каждом шаге $k$ функция $f_k$ может подбираться из рассматриваемого класса функций враждебно по отношению к используемому нами методу генерации последовательности $\{x^k\}$. В частности, $f_k$ может зависеть от

$$\left\{x^1,\xi^1,f_1(\cdot);\ldots;x^{k-1},\xi^{k-1},f_{k-1}(\cdot);x^k\right\}.$$

В стохастической онлайн оптимизации с неточным оракулом можно получить следующие оценки псевдо регрета (см., например, [129], случай неточного оракула в похожем контексте ранее уже частично прорабатывался в [31])

$$\text{(11)} \qquad \min\left\{\mathrm{O}\left(\sqrt{\frac{M^2R^2}{N}} + R\delta\right), \mathrm{O}\left(\frac{M^2\ln N}{\mu N} + R\delta\right)\right\},$$

---

[39] Если $Q = \{x : g(x) \le 0\}$ и на это множество сложно проектироваться, то можно обобщить, сохранив оценку (11), конструкцию прямо-двойственного метода из работы [113] (см. также [44]) на онлайн контекст с таким множеством $Q$.

[40] Это условие можно заменить, считая, что вместо субградиента мы получаем $\delta$-субградиент [2, 131].

где $\|\nabla f_k(x,\xi)\|_* \le M$ – равномерно по $x$, $k$ и п.н. по $\xi$. Эти оценки достигаются (фактически на тех же методах, что и в п. 2 с небольшой оговоркой в сильно выпуклом случае [39, 129]) и неулучшаемы (в том числе для детерминированных постановок с $\delta = 0$ и с линейными функциями $f_k(x)$). Как видно из этих оценок, наличие гладкости не позволяет получить более высокую скорость сходимости. То есть никакого аналога формулы (7) здесь нет. Все что ранее говорилось про прокс-структуру[41] и большие отклонения, насколько нам известно, полностью и практически без изменений переносится и на задачи онлайн оптимизации.

## 4. Стохастические безградиентные и покомпонентные методы с неточным оракулом

Рассматривается задача стохастической выпуклой оптимизации (1)

$$f(x) = E_\xi\left[f(x,\xi)\right] \to \min_{x\in Q}.$$

Предположения те же, что и в первом абзаце п. 2. В частности, п.н. $\|\nabla f(x,\xi)\|_2 \le M$. Здесь важно, что функция $f(x)$ задана не только на множестве $Q$, но и в его $\tau_0$-окрестности (см. ниже), и все предположения делаются не для $x \in Q$, а для $x$ из $\tau_0$-окрестности множества $Q$ (аналогичная оговорка потребуется далее, при перенесении результатов п. 3). Однако теперь оракул не может выдавать стохастический субградиент функции. На каждой итерации мы можем запрашивать у оракула только значения реализации функции $f(x,\xi)$ в нескольких точках $x$. Принципиальная разница есть только между запросом значения (реализации) функции в одной и запросом в двух точках [120, 132]. Здесь мы ограничимся рассмотрением случая двух точек – случай одной точки представляет интерес только в онлайн контексте (см. [131] и цитированную там литературу). Впрочем, есть достаточно большой и популярный класс одноточечных не онлайн постановок, которого мы здесь не касаемся (см., например, [133, 134]).

**Предположение 2.** *$\delta$-оракул выдает (на запрос, в котором указывается только одна точка $x$) $f(x,\xi)+\delta(x,\xi)$, где с.в. $\xi$ независимо разыгрывается из одного и того же распределения, фигурирующего в постановке (1), случайная величина $\delta(x,\xi) = \tilde{\delta}(x) + \bar{\delta}(\xi)$, где $\bar{\delta}(\xi)$ – независимая от $x$ случайная величина с неизвестным распределением (случайность которой может быть обусловлена не только*

[41] За исключением сильно выпуклого случая, для которого нам известны только оценки в евклидовой прокс-структуре. Кроме того, в сильно выпуклом случае в оценках вероятностей больших уклонений $\ln(\ln(N)) \to \ln N$ – доказательство этого утверждения мы не смогли найти (впрочем, см. замечание 4).

*зависимостью от $\xi$), ограниченная по модулю $\delta/2$ (число $\delta$ – допустимый уровень шума), $\tilde{\delta}(x)/(R\delta)$ – неизвестная 1-липшицева функция.*

Далее в изложении мы будем во многом следовать [2, 33, 45, 76, 129–136]. По полученным от оракула зашумленным значениям $f(x,\xi)+\delta(x,\xi)$ мы можем определить стохастический субградиент (важно, что можно обратиться с запросом к оракулу в двух разных точках при одной и той же реализации $\xi$):

$$g_{\tau,\delta}(x,s,\xi)=\frac{n}{\tau}\Big(f(x+\tau s,\xi)+\delta(x+\tau s,\xi)-\big(f(x,\xi)+\delta(x,\xi)\big)\Big)s, \tag{12}$$

где $s$ – случайный вектор (независимый от $\xi$), равномерно распределенный на $S_2^n(1)$ – единичной сфере в 2-норме в пространстве $\mathbb{R}^n$.[42] Из этого представления можно усмотреть, что липшицева составляющая шума $\delta_2$ из предположения 2 и уровень шума $\delta_1$ из предположения 1 связаны соотношением $\delta_2 \sim \delta_1/n$ (см. формулу (24)). В действительности, для обоснования этой связи требуются значительно более громоздкие рассуждения.

Приведем одну из возможных мотивировок введенной в предположении 2 концепции $\delta$-оракула. Предположим, что оракул может считать абсолютно точно значение липшицевой функции, но вынужден нам выдавать лишь конечное (предписанное) число первых бит. Таким образом, в последнем полученном бите есть некоторая неточность (причем мы не знаем по какому правилу оракул формирует этот последний выдаваемый значащий бит). Однако мы всегда можем прибавить (по mod 1) к этому биту случайно приготовленный (независимый) бит. В результате, не ограничивая общности, можно считать, что оракул последний бит выбирает просто случайно в независимости от отброшенного остатка. То что в итоге выдает оракул соответствует концепции $\delta$-оракула.

Перейдем к получению оценок. В отличие от пп. 2, 3 везде далее в этом пункте мы будем считать, что имеет место обратное неравенство на требуемое число итераций $N \geq n$ [132]. Прежде всего, заметим, что[43]

$$E_{s,\xi,\delta}\left[g_{\tau,\delta}(x,s,\xi)\right]=\nabla f_\tau(x)+\nabla_x E_{\tilde{s},\xi,\delta}\left[\delta(x+\tau\tilde{s},\xi)\right], \tag{13}$$

---

[42] С помощью леммы Пуанкаре [138] такой вектор можно сгенерировать за $O(n)$, разыгрывая $n$ независимых одинаково распределенных стандартных нормальных случайных величин и нормируя их на корень из суммы их квадратов.

[43] Взятие математического ожидания по $\delta$ подчеркивает, что $\delta(x,\xi)$ может быть случайной величиной не только потому, что может зависеть от $\xi$, но и потому, что может иметь собственную случайность.

где $\tilde{s}$ – случайный вектор, равномерно распределенный на $B_2^n(1)$ – единичном шаре в 2-норме, а $f_\tau(x) = E_{\tilde{s},\xi}\left[f(x+\tau\tilde{s},\xi)\right]$ – сглаженная[44] версия функции $f(x) = E_\xi\left[f(x,\xi)\right]$. Причем,

(14) $$0 \le f_\tau(x) - f(x) \le M\tau,$$

(15) $$\left\|g_{\tau,\delta}(x,s,\xi)\right\|_2 \le n\left(M + \frac{2\delta}{\tau}\right).$$

Основная идея [2] заключается в подмене задачи (1) следующей задачей

(16) $$f_\tau(x) = E_{\tilde{s},\xi}\left[f(x+\tau\tilde{s},\xi)\right] \to \min_{x\in Q},$$

$\varepsilon/2$-решение которой при $\tau = \varepsilon/(2M)$, будет $\varepsilon$-решением исходной задачи (1).

Считая $\delta = O(\varepsilon)$ и[45] (приведенное условие выполняется, если мы имеем доступ к $\delta$-оракулу из предположения 2)

$$\left\|\nabla_x E_{\tilde{s},\xi,\delta}\left[\delta(x+\tau\tilde{s},\xi)\right]\right\|_2 = O(\varepsilon/R),$$

можно получить для среднего числа итераций (используя те же алгоритмы для задачи (16), что и в п. 2, со стохастическим градиентом (12)), соответствующие аналоги оценок (2), (4):

(17) $$O\left(n^2M^2R^2/\varepsilon^2\right),\ O\left(n^2M^2/(\mu\varepsilon)\right).$$[46]

Если дополнительно известно, что $f(x,\xi)$ – равномерно гладкая по $x$ функция (это условие можно ослабить [130]) и п.н. по $\xi$

---

[44] Все свойства функции $f(x)$ при переходе к $f_\tau(x)$ могут только улучшиться. В частности, $f_\tau(x)$ также выпуклая функция (можно перенести и на сильную выпуклость с не меньшей константой), с константой Липшица и константой Липшица градиента (если таковая существует у $f(x)$) не большей чем у $f(x)$.

[45] Если это условие не выполняется, то все что написано далее останется верным, правда, при более ограничительных условиях на допустимый уровень шума (это касается и всего последующего изложения). Так, если не налагать это ограничение, то потребуется считать $\delta = O\left(\varepsilon^2/\left(\sqrt{n}MR\right)\right)$ или $\delta = O\left(\varepsilon^{3/2}/\left(\sqrt{nL}R\right)\right)$ – в случае, если $f(x)$ имеет $L$-липшицев градиент. Это можно получить с помощью замечания 8.

[46] Аналогично (2), (4) можно переписать оценку (17) не в среднем, как сейчас, а с учетом вероятностей больших уклонений. Это замечание касается и последующих вариаций формулы (17). Нам не известно являются ли оценки (17) оптимальными при заданном уровне шума $\delta = O(\varepsilon)$.

(18) $$\left\|\nabla f(x,\xi)-\nabla f(y,\xi)\right\|_2 \le L\left\|x-y\right\|_2,$$

то вместо (14) будем иметь

(19) $$0 \le f_\tau(x)-f(x) \le \frac{L\tau^2}{2}.$$

Из формулы (19) следует, что можно ослабить требование к неточности: допускать неточность оракула масштаба[47] $\delta \sim M\tau \sim M\sqrt{\varepsilon/L}$.

При сделанных дополнительных предположениях о гладкости (18) за счет ужесточения требований к масштабу допускаемой неточности $\delta$ (как именно требуется это сделать можно усмотреть из формулы (21); ниже мы вернемся к этому вопросу) можно улучшить скорость сходимости (фактор $n^2$ перейдет в $n$):

(20) $$\mathrm{O}\left(nM^2R^2/\varepsilon^2\right),\ \mathrm{O}\left(nM^2/(\mu\varepsilon)\right).$$

Оценки (20) в общем случае не улучшаемы (даже при $\delta = 0$) для гладких стохастических и негладких задач [132]. Фактически это означает, что мы можем выбрать настолько малое $\tau$ (насколько малым мы можем его выбрать определяется $\varepsilon$ и $\delta$), что конечная разность в (12) "превращается" (с нужной точностью) в производную по направлению. Для объяснения отмеченного перехода полезно заметить, что [76, 130–132] (см. также (22))

(21) $$E_{s,\xi}\left[\left\|g_{\tau,\delta}(x,s,\xi)\right\|_2^2\right] \le 4nM^2 + L^2\tau^2n^2 + \frac{8\delta^2n^2}{\tau^2}.$$

**Замечание 7 (техника двойного сглаживания негладких задач Б.Т. Поляка [2], см. также [132, 137]).** За счет подмены изначально негладкого функционала в задаче (1) на

$$f(x) := f_\gamma(x) = E_{\tilde{s}_1,\xi}\left[f(x+\gamma\tilde{s}_1,\xi)\right],\ \gamma \le \varepsilon/(2M),$$

где $\tilde{s}_1$ – случайный вектор (независимый от $\xi$), равномерно распределенный на $B_2^n(1)$, получим новую задачу ($\varepsilon/2$-аппроксимирующую исходную), для которой при достаточно малом $\tau$ будет иметь место оценка (21) с достаточно большим $L \ge 2nM^2/\varepsilon$. Далее решая с помощью уже описанной техники с точностью $\varepsilon := \varepsilon/2$ задачу стохастической оптимизации (1) с

$$\xi := (\tilde{s}_1,\xi),\ f(x,\xi) := f(x+\gamma\tilde{s}_1,\xi),$$

$$g_{\tau,\delta}(x;s_2,\xi) := \frac{n}{\tau}\Big(f(x+\gamma\tilde{s}_1+\tau s_2,\xi)+\delta(x+\gamma\tilde{s}_1+\tau s_2,\xi)-\big(f(x+\gamma\tilde{s}_1,\xi)+\delta(x+\gamma\tilde{s}_1,\xi)\big)\Big)s_2,$$

---

[47] Здесь мы дополнительно считаем, что $\nabla_x E_{\tilde{s},\xi,\delta}\left[\delta(x+\tau\tilde{s},\xi)\right] = 0$. В частности, это условие выполняется, если неточность $\delta(x,\xi)$ имеет независимое от $x$ распределение.

где $s_2$ – случайный вектор (независимый от $\xi$ и $\tilde{s}_1$), равномерно распределенный на $S_2^n(1)$, получим те же оценки (20), только при существенно более жестких условиях на уровень шума $\delta$. К сожалению, получить конструктивное описание этих условий на данный момент не удалось.

Оценки (17) и (20) переносятся и на задачи стохастической онлайн оптимизации (см., например, [43, 129–132, 135]) с возникновением дополнительного фактора $\ln N$ в сильно выпуклом случае (см. формулу (11)). При этом даже в гладком случае не обязательно требовать дополнительно стохастичность исходной постановки для оптимальности оценок (20).

Далее рассматривается не стохастический вариант постановки задачи (1) (не обобщаемый на онлайн постановки) с $Q = \mathbb{R}^n$ (обобщения на произвольные выпуклые множества $Q \subset \mathbb{R}^n$ представляются интересными, но на данный момент нам не известны такие обобщения[48] – в последующих рассуждениях существенным образом используется то, что в точке минимума $\nabla f(x_*) = 0$). Так что теперь $R$ – расстояние от точки старта до решения в 2-норме. В этом варианте, выписанная оценка (21) может быть уточнена

$$(22) \qquad E_s\left[\left\|g_{\tau,\delta}(x,s)\right\|_2^2\right] \le 4n\left\|\nabla f(x)\right\|_2^2 + L^2\tau^2 n^2 + \frac{8\delta^2 n^2}{\tau^2}.$$

Последняя оценка следует из явления концентрации равномерной меры на $S_2^n(1)$ с выделенными полюсами вокруг экватора (см. [138] – в случае покомпонентных методов эта оценка особенно просто получается [43], $s$ в приводимой формуле, и только в ней, соответствует покомпонентной рандомизации):

$$E_s\left[\left\langle \nabla f(x), s\right\rangle^2\right] = \frac{1}{n}\left\|\nabla f(x)\right\|_2^2.$$

Считая для простоты формулировок, что

$$\nabla_x E_{\tilde{s},\delta}\left[\delta\left(x + \tau\tilde{s}\right)\right] = 0,$$

можно распространить метод [87], дающий оценки (8), на текущий контекст, и получить следующие оценки (с $\tau \sim \sqrt{\delta/L}$) числа итераций для достижения точности $\varepsilon$ для случая выпуклой и сильно выпуклой целевой функции соответственно

$$(23) \qquad N_1(\varepsilon) = n \cdot \mathrm{O}\left(\left(\frac{LR^2}{\varepsilon}\right)^{\frac{1}{p+1}}\right), \; N_2(\varepsilon) = n \cdot \mathrm{O}\left(\left(\frac{L}{\mu}\right)^{\frac{1}{p+1}} \ln\left(\frac{LR^2}{\varepsilon}\right)\right)$$

[48] По-видимому, такие обобщения возможны. Также возможно перенесение концепции универсальных методов (см. п. 3) на безградиентные методы и спуски по направлению (покомпонентные спуски) для детерминированных задач (не задач стохастической оптимизации).

при (условия на допустимый уровень шума, при котором оценки (23) имеют такой же вид, с точностью до $\mathrm{O}(1)$, как если бы шума не было[49])

$$\text{(24)} \qquad \delta_1(\varepsilon) \le \frac{1}{n}\mathrm{O}\left(\varepsilon \cdot \left(\frac{\varepsilon}{LR^2}\right)^{\frac{p}{p+1}}\right),\ \delta_2(\varepsilon) \le \frac{1}{n}\mathrm{O}\left(\varepsilon \cdot \left(\frac{\mu}{L}\right)^{\frac{p}{p+1}}\right).$$

По-видимому (строгим доказательством мы не располагаем на данный момент), и в стохастическом случае имеет место аналог формул (23), (24) с заменой в формуле (23)

$$N_1(\varepsilon) = n \cdot \max\left\{\mathrm{O}\left(\left(\frac{LR^2}{\varepsilon}\right)^{\frac{1}{p+1}}\right), \mathrm{O}\left(\frac{DR^2}{\varepsilon^2}\right)\right\},\ N_2(\varepsilon) = n \cdot \max\left\{\mathrm{O}\left(\left(\frac{L}{\mu}\right)^{\frac{1}{p+1}} \ln\left(\frac{LR^2}{\varepsilon}\right)\right), \mathrm{O}\left(\frac{D}{\mu\varepsilon}\right)\right\}.$$

Можно продолжать переносить все написанное в п. 3 на рассматриваемую ситуацию (частично это уже сделано в [43, 120]). Однако мы остановимся лишь на наиболее интересном (на наш взгляд) месте. А именно, на согласовании прокс-структуры с рандомизацией, порождающей сглаживание.

Основным результатом (ввиду замечания 7) в негладком и(или) стохастическом случае с точным оракулом здесь является следующее наблюдение [130, 131]: независимо от выбора прокс-структуры рандомизацию всегда стоит выбирать согласно (12) (если ставить цель – минимизировать число итераций), т.е. с помощью разыгрывания случайного вектора $s$ равномерно распределенного на $S_2^n(1)$. В случае неточного оракула, по-видимому, это утверждение уже перестает быть верным [130]. Ограничимся далее обобщением оценок (20) на случай использования общих прокс-структур.

Приведем соответствующее обобщение формулы (21) (здесь и далее нижний индекс "2" у констант Липшица подчеркивает, что они считаются согласно евклидовой норме из-за сделанного нами выбора способа рандомизации)

$$E_{s,\xi}\left[\left\|g_{\tau,\delta}(x,s,\xi)\right\|_{q'}^2\right] = \mathrm{O}\left(\left(4M_2^2 n + L_2^2\tau^2 n^2 + \frac{8\delta^2 n^2}{\tau^2}\right)E_s\left[\|s\|_{q'}^2\right]\right),$$

где в прямом пространстве выбрана $q$-норма и $1/q + 1/q' = 1$. Согласно замечанию 2 можно считать, что $2 \le q' \le \infty$ – выбирать другие нормы, как правило, бывает не выгодно. Для

[49] В отсутствии шума, оракул нам фактически может выдавать производную по направлению $s$

$$g(x,s) = n\langle\nabla f(x), s\rangle s,$$

точнее $\langle\nabla f(x), s\rangle$, $s$ мы генерируем сами. Если, в свою очередь, считать, что $\langle\nabla f(x), s\rangle$ оракул выдает с аддитивным шумом (для простоты считаем, независящим от $s$) масштаба $\delta := \sqrt{L\delta}$ ($\delta$ в правой части определяется исходя из формулы (24)), то формула (23) останется верной [120].

такого диапазона $E_s\left[\|s\|_{q'}^2\right]=\tilde{O}\left(n^{2/q'-1}\right)$, в частности $E_s\left[\|s\|_{\Omega(\log n)}^2\right]=\tilde{O}\left(n^{-1}\right)$ ($\tilde{O}(\ )$ с точностью до логарифмического фактора от $n$ совпадает с $O(\ )$, аналогично с $\tilde{\Omega}(\ )$). Исходя из такого обобщения, можно привести следующую таблицу, распространяющую оценку (20) на произвольные прокс-структуры ($R^2$ – "расстояние" Брэгмана, согласованное с $q$-нормой, см. замечание 1).

| $f(x)$ – выпуклая | $f(x)$ – $\mu_q$-сильно выпуклая в $q$-норме |
|---|---|
| $O\left(\frac{nM_2^2R^2}{\varepsilon^2}\right)\tilde{O}\left(n^{2/q'-1}\right)$ | $\tilde{O}\left(\frac{nM_2^2}{\mu_q\varepsilon}\right)\tilde{O}\left(n^{2/q'-1}\right)$ |

**Таблица 1**

Выпишем условия, из которых можно получить требования на шум (нам представляется, что здесь это может хорошо прояснить суть дела):

- $\min\left\{M_2\tau, L_2\tau^2/2\right\}=O(\varepsilon)$ – условие достаточной точности аппроксимации исходной функции ее сглаженной версией;
- $L_2^2\tau^2n^2+\frac{8\delta^2n^2}{\tau^2}=O\left(M_2^2n\right)$ – условие "правильной" ограниченности квадрата нормы аппроксимации стохастического градиента.

Выписанные условия позволяют для всех полей таблицы 1 (с оценками) написать соответствующие условия на допустимый уровень шума, и параллельно подобрать оптимальный размер параметра сглаживания $\tau$.

По-видимому, аналогичную таблицу можно записать (на базе конструкций работ [43, 130]) и для оценки (23) ($p=1$)

| $f(x)$ – гладкая выпуклая | $f(x)$ – гладкая $\mu_q$-сильно выпуклая в $q$-норме |
|---|---|
| $\tilde{O}\left(n^{1/q'+1/2}\sqrt{\frac{L_2R^2}{\varepsilon}}\right)$ | $\tilde{O}\left(n^{1/q'+1/2}\sqrt{\frac{L_2}{\mu_q}}\right)$ |

**Таблица 2**

**Замечание 8 (см. [131]).** Если не делать никаких предположений о шуме $\delta(x,\xi)$ в предположении 2 кроме $|\delta(x,\xi)| \le \delta$, то для получения требований на уровень шума $\delta$, потребуется еще воспользоваться следующим утверждением.

Пусть последовательность независимых случайных векторов $\{s_k\}_{k=0}^{N}$, равномерно распределенных на $S_2^n(1)$, и $\left\{x_k\left(\{s_l\}_{l=0}^{k-1}\right)\right\}_{k=1}^{N}$ обладают свойством $E\left[\|x_k - x_*\|_2^2\right] \le R^2$. Тогда

$$E\left[\frac{1}{N}\sum_{k=1}^{N}\left|\langle s_k, x_k - x_*\rangle\right|\right] \le \frac{2R}{\sqrt{n}}.$$

В виду замечаний 2, 4 при использовании этого утверждения можно считать, что $R^2 = \mathrm{O}(V(x_*, x_0))$. Причем константа в $\mathrm{O}(\ )$ может быть сделана $\sim 1$.

Приведенная в первом столбце таблицы 1 оценка при определенных условиях может быть лучше нижней оценки [132].[50] Здесь ситуация аналогична той, о которой написано в конце замечания 2 (будем использовать те же обозначения) и в п. 3 . Например, если $Q = B_1^n(1)$, то в нижней оценке [132] стоит $\tilde{\Omega}\left(nM_1^2/\varepsilon^2\right)$ ($E_\xi\left[\|\nabla f(x,\xi)\|_\infty^2\right] \le M_1^2$), а в таблице будет стоять $\tilde{\mathrm{O}}\left(M_2^2/\varepsilon^2\right)$. Осталось заметить, что $M_2^2 \le nM_1^2$, причем в определенных ситуациях может быть $M_2^2 \ll nM_1^2$.

---

[50] Аналогичное можно сказать и для таблицы 2. А именно, рассмотрим задачу минимизации гладкого выпуклого функционала $f(x)$ с константой Липшица градиента в 2-норме равной $L_2$ на множестве $Q = B_1^n(R)$. Тогда нижняя оценка (при $N \le n$) будет иметь вид (А.С. Немировский, 2015)

$$f(x^N) - f_* \ge \frac{\tilde{C}_1 L_2 R^2}{N^3}.$$

С другой стороны, если использовать обычный FGM с KL-прокс-структурой для этой же задачи, то верхняя оценка будут иметь вид

$$f(x^N) - f_* \le \frac{\tilde{C}_2 L_1 R^2}{N^2},$$

где константа Липшица градиента в 1-норме $L_2/n \le L_1 \le L_2$. Тем не менее, отсюда нельзя сделать вывод, что нижняя оценка достигается. Достигается ли эта нижняя оценка, и, если достигается, на каком методе? – насколько нам известно, это пока открытый вопрос, поставленный А.С. Немировским в 2015 году. Однако если оценивать не число итераций, а общее число арифметических операций и если ограничиться рассмотрением класса функций, для которых стоимость расчета производной по направлению (или значения функции в точке) в $\sim n$ раз меньше стоимости расчета полного градиента [22] (в виду БАД [41, 42] это предположение довольно обременительное, впрочем, если функция задана моделью черного ящика, выдающего только значение функции, а градиент восстанавливается при $n+1$ таком обращении, то сделанное предположение кажется вполне естественным), то при $N \le n$ выписанная выше нижняя оценка (в варианте для общего числа арифметических операций, необходимых для достижения заданной точности) будет соответсвовать оценке из первого столбца таблицы 2.

Основные конкурирующие рандомизации в гладком случае – это рандомизация на евклидовой сфере и покомпонентная рандомизация [43, 130–132], которая используется в основном только с евклидовой прокс-структурой [43]. Исследования последних нескольких лет показали (см., например, [43, 50, 55, 63, 65]), что для довольно большого класса гладких задач выпуклой оптимизации в пространствах огромной размерности,[51] возникающих в самых разных приложениях [50, 55, 60, 63], покомпонентные методы являются наиболее эффективным способом решения (с точки зрения общего числа арифметических операций для достижения заданной точности по функции). Покомпонентные методы, безусловно, заслуживают отдельного подробного обзора. Поэтому здесь мы ограничимся только ссылкой на такие обзоры [43, 139].

Приведем далее несколько примеров, демонстрирующих важность изучения безградиентных методов (часто эти методы называют прямыми методами [2] или методами нулевого порядка [76]).

**Пример 5 (двухуровневая оптимизация [140]).** Требуется решить задачу, возникающую, например, при поиске равновесия по Штакельбергу [106]

$$\psi(x,u) \to \max_{u\in U},$$

$$f(x,u(x)) \to \min_{x}.$$

Из первой задачи находится зависимость $u(x)$, которая входит во вторую (внешнюю) задачу. Проблема здесь в том, что явная зависимость $u(x)$ в общем случае может быть недоступна. Как следствие, могут быть проблемы с расчетом $\nabla u(x)$. Поэтому предлагается приближенно решать первую задачу и использовать безградиентный метод с неточным оракулом для второй. Насколько точно надо решать первую задачу, и какой именно безградиентный метод (с точки зрения чувствительности к неточности) выбирать для второй – определяется сложностью решения первой задачи и свойствами второй.

Рассмотренная двухуровневая задача может быть сильно упрощена, если удается найти ее седловое представление [14, 58, 61, 62, 101, 102, 107–109]. В частности, если функции $\psi(x,u)$, $f(x,u)$ – выпуклы по $x$ и вогнуты по $u$, и $\psi(x,u)$ – простой структуры, то можно (для достаточно большого $\lambda$) заменить исходную задачу на следующую

$$\min_{x}\max_{u\in U}\left[f(x,u)+\lambda\psi(x,u)\right].$$

[51] Во всех этих задачах можно считать полные градиенты, и строить на их базе различные методы. То есть для таких задач выбор покомпонентного метода – осмысленный выбор наиболее быстрого способа решения, а не следствие каких-то (в том числе вычислительных) ограничений на задачи.

Полученную седловую задачу стоит решать методами композитной оптимизации (см. замечание 6 и [9, 14, 101, 102]), чтобы параметр $\lambda$ либо совсем не входил в оценки числа итераций, либо входило очень слабо.

К сожалению, седловое представление возможно далеко не всегда.

**Пример 6 (огромная скрытая размерность).** Пусть $y(x)\in\mathbb{R}^m$, $x\in\mathbb{R}^n$, $n\ll m$. Требуется решить задачу

$$f\left(x,y(x)\right)\to\min_x .$$

Мы предполагаем, что можем эффективно посчитать с необходимой точностью $y(x)$ и $f\left(x,y(x)\right)$ за $\mathrm{O}(m)$. Если для решения этой задачи оптимизации мы будем использовать безградиентный метод, то общее число сделанных арифметических операций пропорционально $mn$ (см. (20), (23)). Заметим, что если бы мы могли использовать обычный градиентный метод, то общее число сделанных арифметических операций также было бы пропорционально $mn$, однако вычисление градиента по разным причинам может быть затруднено (см. пример 5 и метод MCMC [52] для расчета PageRank в подходе [141, 142]). В действительности, часто имеет место следующее полезное наблюдение [141]: если мы можем вычислить значения $y(x)$ и $f\left(x,y(x)\right)$ за $\mathrm{O}(m)$, то мы можем с такой же по порядку сложностью (и затратами памяти) вычислить и производные по фиксированному направлению $h$:

$$\frac{dy(x)}{dx}h,\ \left\langle\frac{\partial f(x,y)}{\partial x},h\right\rangle+\left\langle\frac{\partial f(x,y)}{\partial y},\frac{dy(x)}{dx}h\right\rangle.$$

Тем не менее, тут требуется много оговорок, в том числе про точность расчетов. Если не вдаваться в детали, то такие рассуждения также приводят к затратам пропорциональным $mn$, где $n$ возникло в виду оценок (20), (23) для спусков по направлению. Только в отличие от полно-градиентного метода для покомпонентного метода константа Липшица градиента функционала в оценке числа итераций уже будет рассчитываться не по худшему направлению, а в среднем (это может давать выгоду, по порядку равную корню квадратному из размерности пространства [43, 143]), да и, как правило, будет ощутимая выгода в затрачиваемых ресурсах машинной памяти [141]. Оговорки о точности здесь все же необходимы, поскольку для безградиентных методов и спусков по направлению требования к точности могут существенно отличаться (об этом ранее уже было немного написано в данном пункте). Как следствие, в оценку $\mathrm{O}(m)$ необходимо явно вводить зависимость от точности вычисления $y(x)$ и $f\left(x,y(x)\right)$. Об этом планируется написать отдельно.

В примерах 5, 6, в действительности, требуются некоторые оговорки о невозможности или неэффективности использования БАД для полно-градиентного метода (см. п. 2, а также [41, 42]). Нам известны случаи, подпадающие под разобранные примеры, в которых не понятно, как можно было бы воспользоваться БАД [106, 141, 142]. В частности, в работах [141, 142], соответствующей примеру 6, сложность в том, что БАД хочется использовать для ускорения вычисления матрицы Якоби отображения (вектор функции) $y(x) \in S_m(1)$, неявно заданного уравнением $y = P(x)y$, со стохастической (по столбцам) эргодической матрицей $P(x)$ со спектральной щелью $\alpha$ и числом ненулевых элементов $sm$. Метод простой итерации позволяет с точностью $\varepsilon$ найти $y(x)$ за время $\mathrm{O}\left(sm\left(n+\alpha^{-1}\ln\left(\varepsilon^{-1}\right)\right)\right)$ с затратами памяти $\mathrm{O}(sm)$. Нам не известно более эффективного способа расчета матрицы Якоби отображения $y(x)$, чем естественное обобщение метода простой итерации (для продифференцированного по $x$ уравнения $y(x) = P(x)y(x)$), требующее затрат времени $\mathrm{O}\left(smn\alpha^{-1}\ln\left(\varepsilon^{-1}\right)\right)$ и памяти $\mathrm{O}(smn)$. Для реальных приложений [141, 142]: $m \sim 10^9$, $s \sim 10^2$, $n \sim 10^3$, $\alpha \sim 10^{-1}$, $\varepsilon \sim 10^{-12}$ Отсюда ясно, что при использовании полно-градиентного метода просто невозможно будет выделить даже у 64-битной операционной системы, стоящей на самом современном персональном компьютере, необходимой памяти под работающую программу, в основе которой лежит полно-градиентный подход.

Кроме того, в примерах 5, 6 важно уметь эффективно пересчитывать значения $u(x)$ или $y(x)$, а не рассчитывать их каждый раз заново (на каждой итерации внешнего цикла). Поясним сказанное. Предположим, что мы уже как-то посчитали, скажем, $u(x)$, решив с какой-то точностью соответствующую задачу оптимизации. Тогда для вычисления $u(x+\Delta x)$ (на следующей итерации внешнего цикла) у нас будет хорошее начальное приближение $u(x)$. А как известно (см. пп. 2–4) расстояние от точки старта до решения (не в сильно выпуклом случае) существенным образом определяет время работы алгоритма оптимизции. Эта конструкция (hot/warm start) напоминает фрагмент обоснования сходимости методов внутренней точки при изучении движения по центральному пути [8–11]. Тем не менее, известные нам приложения (пример 4 п. 3 и [61, 62, 141, 142]) пока как раз всецело соответствуют сильно выпуклой ситуации. Связано это с тем, что если расчет $u(x)$ или $y(x)$ с точностью $\varepsilon$ осуществляется за $\mathrm{O}\left(\mathrm{C}\ln\left(R/\varepsilon\right)\right)$ операций, то для внешней задачи можно выбирать самый быстрый метод (а, стало быть, и самый требовательный к точности), и с точностью до того, что стоит под логарифмом, общая трудоемкость будет просто прямым произведением трудоемкостей решения внутренней и внешних задач по отдельности. Как правило, такое сочетание оказывается недоминируемым.

Также необходимо отметить, что, как правило, итоговые задачи оптимизации (после подстановки зависимости $u(x)$ или $y(x)$ в задачу верхнего уровня) в этих примерах получаются не выпуклыми. В этой связи можно лишь говорить о локальной сходимости к стационарной точке. В отсутствие выпуклости даже если ограничиться локальной сходимостью многое из того, что описано в данной статье, требует отдельного рассмотрения [76, 96].

Отметим в заключение, что если немного по-другому посмотреть на описанное в этом пункте, то можно заметить следующее. Какой бы большой (но равномерно ограниченный по итерациям) шум ни был, если $\delta\left(x+\tau s,\xi\right)$ имеет распределение не зависящее от $s$, то (возможно через очень большое число итераций) мы сможем сколь угодно точно (по функции) решить задачу! Аналогичное можно сказать, если мы изначально исходим из концепции оракула, выдающего зашумленное значение $\left\langle \nabla f\left(x\right),s\right\rangle$, причем зашумленность не зависит от $s$. Все это восходит к идеям Р. Фишера, развитым О.Н. Граничиным и Б.Т. Поляком [33].

## 5. Заключение

Авторы выражают благодарность А.С. Антипину, С. Бойду, А.Ю. Горнову, Ф. Глинёру, О.Н. Граничину, А.И. Голикову, О. Деволдеру, Н.К. Животовскому, А.В. Назину, Е.А. Нурминскому, А. Рахлину, П. Рихтарику, А. Содину, В.Г. Спокойному, С.П. Тарасову, Э. Хазану, С.В. Чуканову, А.Г. Шапиро, А.Б. Юдицкому и особенно А.С. Немировскому за возможность обсуждения разных частей данного текста в разное время, а также Ильнуре Усмановой, обнаружившей несколько опечаток.

Особо хотелось бы отметить то позитивное влияние, которое на каждого из авторов оказал Борис Теодорович Поляк, отметившим в 2015 году свое восьмидесятилетие [144]. Так Ю.Е. Нестеров защищал кандидатскую диссертацию под руководством Б.Т. Поляка в 1983 г. Два других соавтора в значительной степени являются учениками Ю.Е. Нестерова по части, связанной с численными методами выпуклой оптимизации. Набор идей, здоровый оптимизм, понятность и открытость Бориса Теодоровича всегда делает общение и работу с ним чрезвычайно интересной и полезной. Настоящая статья во многом была вдохновлена его книгой [2], и попыткой ее переосмыслить с современных позиций. Несмотря на то, что большей частью эта книга была написана почти 40 лет назад, она до сих пор не утратила свою актуальность, о чем говорит, например, недавний выход ее нового издания.

Настоящая статья, в основу которой положено выступление [145], была написана уже к осени 2014 года. К сожалению, из-за объема и стиля у нас возникли проблемы с ее публикацией. Мы хотели бы поблагодарить Леонида Васильевича Стрыгина любезно согласившегося взять данную статью (несмотря на ее большой объем) в журнал Труды МФТИ в марте 2016 года.

Исследование А.В. Гасникова и П.Е. Двуреченского в части 3 выполнено в ИППИ РАН за счет гранта Российского научного фонда (проект №14-50-00150), а в части 2, 4 при поддержке гранта РФФИ 15-31-20571-мол_а_вед; исследование Ю.Е. Нестерова выполнено при поддержке гранта РФФИ 14-01-00722-а.

## СПИСОК ЛИТЕРАТУРЫ